\documentclass[10pt]{article}
\usepackage{mystyle}
\usepackage{subfiles}

\begin{document}

\title{A Divergence-Conforming Hybridized Discontinuous Galerkin Method for the Incompressible Reynolds Averaged Navier-Stokes Equations}
\author{Eric L. Peters and John A. Evans}
\date{}
\maketitle

\section*{Abstract}

\noindent We introduce a hybridized discontinuous Galerkin method for the incompressible Reynolds Averaged Navier-Stokes equations coupled with the Spalart-Allmaras one equation turbulence model.  With a special choice of velocity and pressure spaces for both element and trace degrees of freedom, we arrive at a method which returns point-wise divergence-free mean velocity fields and properly balances momentum and energy.  We further examine the use of different polynomial degrees and meshes to see how the order of the scalar eddy viscosity affects the convergence of the mean velocity and pressure fields, specifically for the method of manufactured solutions.  As is standard with hybridized discontinuous Galerkin methods, static condensation can be employed to remove the element degrees of freedom and thus dramatically reduce the global number of degrees of freedom.  Numerical results illustrate the effectiveness of the proposed methodology.

\noindent \textbf{Key Words:} Hybridized discontinuous Galerkin methods; Divergence conforming discretizations; RANS: Reynolds averaged Navier-Stokes; Incompressible flow; Spalart-Allmaras turbulence model; Static condensation

\section{Introduction}
The discontinuous Galerkin (DG) method was originally introduced by Reed and Hill in \cite{Reed_1973} to solve the neutron transport equation, however with little analysis of its properties.  A more rigorous exploration of the method and its properties were discovered and discussed in \cite{ johnson_analysis_nodate,lasaint_finite_1974}.  Since then, the method has gained popularity and has since been extended to a multitude of other problems governed by partial differential equations (PDEs) detailed in \cite{cockburn_book} and \cite{hesthaven_nodal_2008}.  DG methods combine the advantages of classical finite volume methods (FVM) with finite element methods (FEM).  Much like FVM, DG has a natural framework for dealing with advection problems, or more generally hyperbolic problems, but in addition it has a natural extension to higher order methods, in a similar fashion to FEM.  Thus DG is capable of resolving solutions with large gradients, including shocks, as well as dealing with complex geometries.  Another important advantage of having a method which is higher order is that more work can be done on processor for element interiors, limiting the overall percentage of work needed for communication in a parallel framework.  In fact, it has been regularly used to solve large-scale forward \cite{breuer_sustained_2014,wilcox_high-order_2010} and inverse problems \cite{bui-thanh_extreme-scale_2012}.  However, there is a major drawback associated with this class of method.  In general, DG methods introduce an increased amount of degrees of freedom (DOFs) which cannot be alleviated by increasing the polynomial order.  Consequently, implicit time integration schemes are often intractable and one must resort to explicit schemes hindering the range of acceptable Courant-Friedrichs-Lewy (CFL) condition numbers.  This can make the method prohibitively expensive in comparison to other existing numerical methods, such as continuous Galerkin (CG) or spectral element methods.\\
\\
Within the last decade, Cockburn and his collaborators have devised a method which alleviates the high cost of DG methods using hybridization.  The resulting approach, referred to as the hybridized discontinuous Galerkin (HDG) method, has been applied successfully to several types of PDEs including, but not limited to, Poisson's equation \cite{cockburn_unified_2009}, advection-diffusion equations \cite{cockburn_hybridizable_2009,egger_hybrid_2010,nguyen_implicit_2009,nguyen_implicit_2009-1}, the Stokes equations \cite{cockburn_derivation_2009,cockburn_comparison_2010,cockburn_analysis_2011,nguyen_hybridizable_2010}, and even the Euler and Navier-Stokes equations \cite{moro_navier-stokes_nodate,peraire_hybridizable_2010}.  With an HDG method, there are both interior DOFs that reside on the interiors of elements and trace (or skeleton) DOFs that reside on the edges of polygons in 2-D and faces of polyhedra in 3-D.  The trace DOFs act as communicators across element boundaries.  In fact, static condensation can be employed to write the interior DOFs in terms of the trace DOFs, allowing one to remove the interior DOFs entirely from the discrete system of equations.  In this manner, the HDG method allows one to significantly reduce the number of DOFs in the global system as compared with a corresponding DG method, rendering the HDG method competitive with higher-order CG and spectral element methods.  Once the trace DOFs are determined, the interior DOFs can be realized locally in an element-by-element manner.  Furthermore, with the supposition of convergence of a DG method, the equivalent HDG method is guaranteed to converge \cite{cockburn_projection-based_2010}, with this condition being sufficient but not necessary.\\
\\
\noindent Recently, there has been a concerted effort in further developing HDG methods for the incompressible Navier-Stokes equations and related PDE's, including the development of stable mixed elements and structure-preserving discretizations \cite{rhebergen_hybridizable_2017} as well as super convergent post-processing procedures \cite{nguyen_implicit_2011}.
There has additionally been efforts to alleviate the associated complexity of developing new HDG discretizations.
Bui-Thanh has contributed to this avenue by constructing HDG methods from a Godunov approach \cite{bui-thanh_godunov_2015} as well as from the Rankine-Hugoniot condition \cite{bui-thanh_rankine-hugoniot_2015}.  These types of methodologies render the trace unknowns as upwind states, while providing a parameter-free framework.\\
\\
In this work, we construct a new HDG method for the incompressible Reynolds-Averaged Navier-Stokes (RANS) equations.  To construct our method, we begin with the divergence-conforming HDG method recently introduced for the incompressible Navier-Stokes equations by Rhebergen and Wells in \cite{rhebergen_hybridizable_2017}, and we modify the method to handle a RANS closure model.  More specifically, we design an HDG method for the Spalart-Allmaras one-equation turbulence model \cite{spalart_one-equation_1992}, where we treat the transport equation for the eddy viscosity as an advection-diffusion equation.  It is important to note that we introduce trace DOFs for the primal variables, that is, the flow field and turbulence variables, for the sake of implementational ease.  This is in contrast with many alternative HDG methods where trace DOFs represent inter-element fluxes.\\
\\
A defining feature of our HDG method is that it produces a point-wise divergence-free mean velocity field.  It has been shown in recent works that discretizations of the incompressible Navier-Stokes equations which exactly preserve the divergence-free constraint are typically more robust than methods which only satisfy the divergence-free constraint in a weak manner \cite{divfreereview}.  For example, divergence-conforming discretizations simultaneously conserve momentum and energy \cite{evans_iso1}, while discretizations which only weakly satisfy the divergence-free constraint typically conserve either momentum or energy but not both.  The velocity error is also decoupled from the pressure error for divergence-conforming discretizations, a property commonly referred to as pressure robustness \cite{dgdivfree}.  Finally, mass conservation is considered to be of prime importance for coupled-flow transport \cite{stabadvdiff}, so divergence-free discretizations are particularly attractive for such applications.  This last point inspires our current work, as the Spalart-Allmaras model may be viewed as a coupled-flow transport problem wherein the eddy viscosity is a scalar transported by the fluid flow.

An outline of this paper is as follows.  In the next section, we recall the strong form of the incompressible RANS equations.  In Section 3, we present a semi-discrete divergence-conforming HDG method for the incompressible RANS equations, assuming an underlying linear eddy viscosity model, and we show that the method conserves mass in a point-wise manner, conserves momentum globally, and is energy stable.  Motivated by the fact that most turbulence model equations are transport equations, we recall the strong form of the scalar advection-diffusion equation in Section 4, and we present a semi-discrete HDG method for the advection-diffusion equation, assuming an underlying divergence-free velocity field, in Section 5.  In Section 6, we present a semi-discrete HDG method for the incompressible RANS equations with the Spalart-Allmaras one-equation turbulence model by combining our previous semi-discrete HDG methods for the incompressible RANS equations and the advection-diffusion equation.  In Section 7, we present a staggered time-integration scheme for our HDG method for the incompressible RANS equations with the Spalart-Allmaras model, and we illustrate how static condensation can be employed to reduce computational cost in Section 8.  We present illustrative numerical results in Section 9, and finally, in Section 10, we draw conclusions and discuss future research directions.

\section{Strong Form of the Reynolds Averaged Navier-Stokes Equations}

We begin this paper by recalling the strong form of the incompressible RANS equations.  The derivation of the incompressible RANS equations relies on a \textit{Reynolds decomposition} of the velocity and pressure fields into \textit{mean} and \textit{fluctuating components}, viz.:
\begin{equation}
\begin{aligned}
\bf{u} &= \overline{\bf{u}} +\bf{u}'\\
p &=  \overline{p}  + p'
\end{aligned} 
\label{mean}
\end{equation}
where $\overline{*}$ and $*'$ denote the mean and fluctuation respectively.  We then obtain the incompressible RANS equations by taking the mean of the incompressible Navier-Stokes equations and exploiting the above Reynolds decompositions.  The resulting system is displayed below:

\begin{equation}
\begin{aligned}
& \frac{\partial \overline{\bf{u}}}{\partial t} +  \nabla \cdot (\overline{\bf{u}} \otimes \overline{\bf{u}}) + \frac{1}{\rho}  \nabla \overline{p} -  \nabla \cdot (2 \nu \nabla^S \bf{u}) + \nabla \cdot (\overline{\bf{u}' \otimes \bf{u}'}) = \bf{f} \\
& \nabla \cdot \overline{\bf{u}} = 0
\end{aligned} 
\label{NS}
\end{equation}
Above, $\nabla^S \bf{u} = \frac{1}{2}(\nabla \overline{\bf{u}} + (\nabla \overline{\bf{u}})^T)$ denotes the \textit{mean strain rate tensor} and $\overline{\bf{u}' \otimes \bf{u}'}$ denotes the \textit{Reynolds stress tensor}.  

In this paper, we will focus our attention on \textit{linear eddy viscosity models}.  These models assume that the mean strain rate tensor and the anisotropic part of the Reynolds stress tensors are aligned, and hence there exists an \textit{eddy viscosity} $\nu_T$ such that:
\begin{equation}
\begin{aligned}
\overline{\bf{u}' \otimes \bf{u}'} = -2 \nu_T \nabla^S {\overline{\bf{u}}} + \frac{2}{3} k \bf{I}
\end{aligned}
\end{equation} 
where $k = \frac{1}{2} \textup{tr}\left( \overline{\bf{u}' \otimes \bf{u}'} \right)$ is the \textit{turbulent kinetic energy} and $\bf{I}$ is the identity tensor.  It is generally assumed that the eddy viscosity is non-negative.  Otherwise, the resulting system may be unstable.  In what follows, we combine the contributions of the mean pressure field and the isotropic part of the Reynolds stress tensor using a \textit{modified pressure} $\overline{p} + \frac{2}{3} \rho k$.  With an abuse of notation, we use $\overline{p}$ to denote this modified pressure.

We are now ready to state the strong form of the incompressible RANS equations subject to a linear eddy viscosity model.  Let $\Omega$ denote a Lipschitz and bounded domain in $\mathbb{R}^d$, where $d=2$ for two-dimensional domains and $d=3$ for three-dimensional domains.  Let $\partial \Omega$ denote the boundary of $\Omega$ with outward unit normal ${\bf{n}}$.  Let $\partial \Omega$ be partitioned into a Dirichlet boundary $\Gamma_D$ and a Neumann boundary $\Gamma_N$ such that $\partial \Omega = \overline{\Gamma_D \cup \Gamma_N}$ and $\Gamma_D \cap \Gamma_N = \emptyset$.  We assume a constant density $\rho \in \mathbb{R}^+$, a variable kinematic viscosity $\nu: \Omega \times (0,\infty) \rightarrow \mathbb{R}^+$, a variable body force ${\bf{f}}: \Omega \times (0,\infty) \rightarrow \mathbb{R}^d$, a variable velocity specification on the Dirichlet boundary ${\bf{g}}: \Gamma_D \times (0,\infty) \rightarrow \mathbb{R}^d$, a variable traction specification on the Neumann boundary ${\bf{h}}: \Gamma_N \times (0,\infty) \rightarrow \mathbb{R}^d$, and a variable initial velocity $\overline{\bf{u}}_{0}: \Omega \rightarrow \mathbb{R}^d$.  We also \textit{freeze} the turbulent viscosity in the following presentation.  That is, we assume that $\nu_T: \Omega \times (0,\infty) \rightarrow \mathbb{R}^+$ is a known function.  This will allow us to examine the properties of our discretization for the incompressible RANS equations independent of the chosen eddy viscosity model.  We will later make the turbulent viscosity a function of the unknown mean flow field and unknown turbulence variables.  With the above assumptions, the strong form is as follows:
\begin{mdframed}
\textbf{Strong Form  for the Incompressible RANS Equations}\\
\\
Find $\overline{\bf{u}} : \bar{\Omega} \times [0,\infty) \rightarrow \mathbb{R}^d$ and $\overline{p}: \Omega \times (0,\infty) \rightarrow \mathbb{R}$ such that:
\begin{equation}
\begin{aligned}
\frac{\partial \overline{\bf{u}}}{\partial t} +  \nabla \cdot (\overline{\bf{u}} \otimes \overline{\bf{u}}) + \frac{1}{\rho}  \nabla \overline{p} -  \nabla \cdot (2( \nu+\nu_T) \nabla^S \overline{\bf{u}}) & = \bf{f} && \text{ in } \Omega \times (0,\infty) \\
\nabla \cdot \overline{\bf{u}} & = 0 && \text{ in } \Omega \times (0,\infty)\\
\overline{\bf{u}} & = {\bf{g}} && \text{ on } \Gamma_D \times (0,\infty)\\
\left[-\frac{1}{\rho}\overline{p}\textbf{I} +  2 (\nu+\nu_T) \nabla^S \overline{\bf{u}}\right] \cdot {\bf{n}}  - \text{min}(\overline{\bf{u}} \cdot {\bf{n}},0) \overline{\bf{u}} & = {\bf{h}} && \text{ on } \Gamma_N \times (0,\infty)\\
\overline{\bf{u}}(\cdot,0) & = \overline{\bf{u}}_{0}  && \text{ in } \Omega\\
\end{aligned}\label{NS_S}
\end{equation}
\end{mdframed}

Note that we have utilized somewhat unorthodox Neumann boundary conditions above.  On the outflow parts $\Gamma_N$ (i.e., where $\overline{\bf{u}} \cdot {\bf{n}}  \geq 0$), we set the traction, i.e., $\left[-\frac{1}{\rho}\overline{p}\textbf{I} +  2 (\nu+\nu_T) \nabla^S \overline{\bf{u}}\right] \cdot {\bf{n}} = {\bf{h}}$, as is standard.  However, on the inflow parts of $\Gamma_N$, we instead set the sum of momentum, pressure, and diffusive fluxes, i.e., $\left[-(\overline{\bf{u}} \otimes \overline{\bf{u}}) - \frac{1}{\rho}\overline{p} \textbf{I} +  2 (\nu+\nu_T) \nabla^S \overline{\bf{u}}\right] \cdot  {\bf{n}} = {\bf{h}}$.  This yields a well-posed formulation in the presence of backflow \cite{donothing}.


\section{Semi-Discrete HDG Formulation for the Reynolds Averaged Navier-Stokes Equations}
\label{sec:semidiscrete_RANS}
\begin{figure}[b!]
\hspace*{0cm}
  \centering
  \includegraphics[scale=1.15]{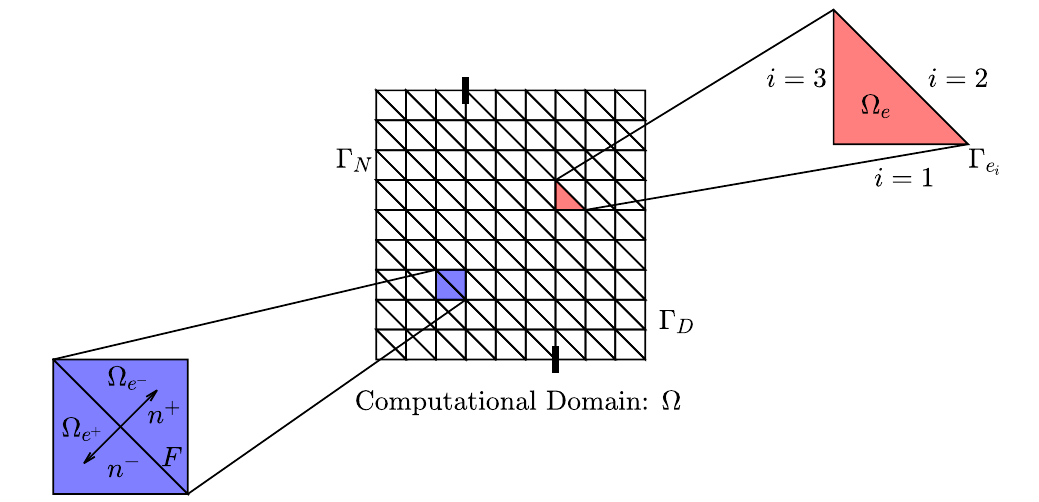}
\caption{Mesh Notation.}
\label{notation}
\end{figure}

We are now ready to construct an HDG method for the incompressible RANS equations subject to a linear eddy viscosity model.  In this section, we discretize in space, and later, we discretize in time.  To discretize in space, we first introduce a \textit{mesh} over which the incompressible RANS equations will be discretized.  Let $\mathcal{T} := \left\{ \Omega_e \right\}_{e=1}^{\text{nel}}$ be a triangulation of the domain into non-overlapping simplex elements such that $\Omega = \overline{ \cup_{e=1}^{\text{nel}} \Omega_e}$.  The $i^{\text{th}}$ edge or face of an element $\Omega_e \subset \mathcal{T}$ is denoted by $\Gamma_{e_i}$ and the outward unit normal vector on $\Gamma_{e_i}$ to $\Omega_{e}$ is denoted by ${\bf{n}}$.  Two adjacent cells $\Omega_{e^+}$ and $\Omega_{e^-}$ share an edge or face $F$, and we call such an $F$ an \textit{interior facet}.  Any edge or face $F$ that lies on the boundary of the domain $\partial \Omega$ is called a \textit{boundary facet}.  The sets of interior and boundary facets are denoted by $\mathcal{F}_{\text{int}}$ and $\mathcal{F}_{\text{bdy}}$ respectively.  The set of all facets is denoted by $\mathcal{F} := \mathcal{F}_{\text{int}} \cup \mathcal{F}_{\text{bdy}}$, and we define the \textit{mesh skeleton} as $\tilde{\Gamma} = \cup_{F \in \mathcal{F}} F$.  A pictorial representation of the above objects can be seen in Fig. \ref{notation}.
\\
\\
Let $P_{k}(D)$ denote the space of polynomials of degree $k \geq 0$ on a domain $D$.  For polynomial degree $k \geq 1$, we consider the following finite element spaces defined on the mesh:
\begin{equation}
\begin{aligned}
V^h &:= \{ {\bf{v}}^h \in [L^2(\Omega)]^d: \left.{\bf{v}}^h\right|_{\Omega_e} \in [P_k(\Omega_e)]^d \hspace{2mm} \forall  \Omega_e \in \mathcal{T} \}\\
\hat{V}^h &:= \{ \hat{\bf{v}}^h \in [L^2(\tilde{\Gamma})]^d: \left.\hat{\bf{v}}^h\right|_{F} \in [P_{k}(F)]^d \hspace{2mm} \forall F \in \mathcal{F} \}\\
Q^h &:= \{q^h \in L^2(\Omega): \left.q^h\right|_{\Omega_e} \in P_{k-1}(\Omega_e)  \hspace{2mm}  \forall  \Omega_e \in \mathcal{T} \} \\
\hat{Q}^h &:= \{ \hat{q}^h \in L^2(\tilde{\Gamma}): \left.\hat{q}^h\right|_{F} \in P_{k}(F) \hspace{2mm} \forall F \in \mathcal{F}\}
\end{aligned}
\label{spaces}
\end{equation}
\\
We will approximate the velocity and pressure fields over element interiors using the spaces $V^h$ and $Q^h$, respectively.  We will approximate the velocity and pressure fields over the mesh skeleton using the spaces $\hat{V}^h$ and $\hat{Q}^h$, respectively.  We will henceforth refer to $V^h$ and $Q^h$ as the (discrete) velocity and pressure spaces and we will refer to $\hat{V}^h$ and $\hat{Q}^h$ as the (discrete) velocity and pressure \textit{trace} spaces.  It is important to note that functions in the velocity and pressure spaces are defined on the whole domain $\Omega$, whereas functions in the trace spaces are defined only on the mesh skeleton $\tilde{\Gamma}$.  We also introduce the spaces:
\begin{equation}
\begin{aligned}
\hat{V}_g^h &:= \{ \hat{\bf{v}}^h \in \hat{V}^h:  \hat{\bf{v}}^h = \textbf{g} \text{ on } \Gamma_D \}\\
\hat{V}_0^h &:= \{ \hat{\bf{v}}^h \in \hat{V}^h:  \hat{\bf{v}}^h = \textbf{0} \text{ on } \Gamma_D \}
\end{aligned}
\label{spaces}
\end{equation}
which will correspond to our velocity trace trial and test spaces.\\
\\
Note that functions in the velocity and pressure spaces $V^h$ and $Q^h$ may be discontinuous across cell boundaries.  Let $F$ be an interior facet of the mesh shared by two elements $\Omega_{e^+}$ and $\Omega_{e^-}$.  We denote the outward facing normals on $F$ to $\Omega_{e^+}$ and $\Omega_{e^-}$ as $\bf{n}^+$ and $\bf{n}^-$, respectively.  For ${\bf{y}} \in V^h$, we denote the trace of ${\bf{y}}|_{\Omega_{e^+}}$ along $F$ as ${\bf{y}}^+$, and we denote the trace of ${\bf{y}}|_{\Omega_{e^-}}$ along $F$ as ${\bf{y}}^-$.  The \textit{jump operator} across $F$ can then be defined as $\llbracket {\bf{y}} \rrbracket = {\bf{y}}^+ \cdot {\bf{n}}^+ + {\bf{y}}^- \cdot {\bf{n}}^-$.\\
\\
With all of the aforementioned terminology in place, we are now ready to state our semi-discrete HDG formulation for the incompressible RANS equations.  Our method is the natural extension of the divergence-conforming HDG method of Rhebergen and Wells \cite{rhebergen_hybridizable_2017} to the RANS setting, wherein advective fluxes are treated using upwinding and diffusive fluxes are treated using the symmetric interior penalty method.\\
\begin{mdframed}
\textbf{Semi-Discrete HDG Formulation for the Incompressible RANS Equations}\\
\\
Find $\left( \overline{{\bf{u}}}^h, \hat{\overline{{\bf{u}}}}^h, \overline{p}^h,  \hat{\overline{p}}^h \right) \in V^h \times \hat{V}_g^h \times Q^h \times \hat{Q}^h$ such that:\\
\\
\noindent \textit{Continuity Equation}
\begin{equation}
\begin{aligned}
\sum_e \int_{\Omega_e}\frac{1}{\rho}  \nabla \cdot \overline{{\bf{u}}}^h q^h d\Omega = 0 \hspace{80mm} \forall q^h \in Q^h
\end{aligned}
\label{mass}
\end{equation}
\noindent \textit{Continuity Conservativity Condition}
\begin{equation}
\begin{aligned}
\sum_e \sum_i \int_{\Gamma_{e_i}}\frac{1}{\rho} \overline{{\bf{u}}}^h \cdot {\bf{n}} \hat{q}^h d\Gamma - \int_{\partial \Omega}\frac{1}{\rho} \hat{\overline{{\bf{u}}}}^h \cdot {\bf{n}} \hat{q}^h d\Gamma = 0 \hspace{44.5mm} \forall \hat{q}^h \in \hat{Q}^h
\end{aligned}
\label{masscon}
\end{equation}
\noindent \textit{Momentum Equation}
\begin{equation}
\begin{aligned}
&\int_{\Omega} \frac{\partial \overline{{\bf{u}}}^h}{\partial t} \cdot {\bf{v}}^h d \Omega \\
-\sum_e &\int_{\Omega_e} (\overline{{\bf{u}}}^h \otimes \overline{{\bf{u}}}^h) : \nabla {\bf{v}}^h d \Omega
- \sum_e \int_{\Omega_e} \frac{1}{\rho} \overline{p}^h \textbf{I} :  \nabla {\bf{v}}^h d \Omega
+ \sum_e \int_{\Omega_e}  2 (\nu +\nu_T) \nabla^S \overline{{\bf{u}}}^h  :  \nabla^S {{\bf{v}}}^hd \Omega \\
+\sum_e \sum_i &\int_{\Gamma_{e_i}} (\overline{{\bf{u}}}^h \otimes \overline{{\bf{u}}}^h) : ({\bf{v}}^h \otimes {\bf{n}}) d\Gamma
+\sum_e \sum_i \int_{\Gamma_{e_i}} (\hat{\overline{{\bf{u}}}}^h - \overline{{\bf{u}}}^h) \otimes \lambda \overline{{\bf{u}}}^h : ({\bf{v}}^h \otimes {\bf{n}}) d\Gamma\\
+\sum_e \sum_i &\int_{\Gamma_{e_i}}\frac{1}{\rho} \hat{\overline{p}}^h \textbf{I} : ({\bf{v}}^h \otimes {\bf{n}}) d\Gamma
-\sum_e \sum_i \int_{\Gamma_{e_i}}2 (\nu +\nu_T) \nabla^S \overline{{\bf{u}}}^h: ({\bf{v}}^h \otimes {\bf{n}}) d\Gamma \\
 -\sum_e \sum_i &\int_{\Gamma_{e_i}}{\frac{2C_{pen}}{h_e} (\nu + \nu_T)} (\hat{\overline{{\bf{u}}}}^h - \overline{{\bf{u}}}^h)   \otimes {\bf{n}}: ({\bf{v}}^h \otimes {\bf{n}}) d\Gamma \\
 + \sum_e \sum_i &\int_{\Gamma_{e_i}}2 (\nu+\nu_T)[ (\hat{\overline{{\bf{u}}}}^h - \overline{{\bf{u}}}^h)  \otimes {\bf{n}}] :  \nabla^S {{\bf{v}}}^h d\Gamma
 - \int_{\Omega} {\bf{f}} \cdot {\bf{v}}^h d\Omega = 0 \hspace{20mm}
\forall {\bf{v}}^h \in {V}^h
\end{aligned}
\label{mom}
\end{equation}
\noindent \textit{Momentum Conservativity Condtion}
\begin{equation}
\begin{aligned}
-\sum_e \sum_i &\int_{\Gamma_{e_i}} (\overline{{\bf{u}}}^h \otimes \overline{{\bf{u}}}^h) : (\hat{{{\bf{v}}}}^h \otimes {\bf{n}}) d\Gamma
-\sum_e \sum_i \int_{\Gamma_{e_i}} (\hat{\overline{{\bf{u}}}}^h - \overline{{\bf{u}}}^h)  \otimes \lambda \overline{{\bf{u}}}^h : (\hat{{{\bf{v}}}}^h \otimes {\bf{n}}) d\Gamma \\
-\sum_e \sum_i &\int_{\Gamma_{e_i}}\frac{1}{\rho} \hat{\overline{p}}^h \textbf{I} : (\hat{{{\bf{v}}}}^h \otimes {\bf{n}}) d\Gamma
+\sum_e \sum_i \int_{\Gamma_{e_i}}2 (\nu +\nu_T) \nabla^S\overline{{\bf{u}}}^h:  (\hat{{{\bf{v}}}}^h \otimes {\bf{n}}) d\Gamma \\
+\sum_e \sum_i &\int_{\Gamma_{e_i}}{\frac{2C_{pen}}{h_e} (\nu + \nu_T)} (\hat{\overline{{\bf{u}}}}^h - \overline{{\bf{u}}}^h)   \otimes {\bf{n}}:  (\hat{{{\bf{v}}}}^h \otimes {\bf{n}}) d\Gamma \\
+&\int_{\Gamma_N} (1-\lambda) (\hat{\overline{{\bf{u}}}}^h \cdot {\bf{n}}) (\hat{\overline{{\bf{u}}}}^h \cdot \hat{{{\bf{v}}}}^h) d \Gamma
+ \int_{\Gamma_N} {\bf{h}} \cdot \hat{{{\bf{v}}}}^h d \Gamma = 0
\hspace{28.5mm }\forall \hat{{{\bf{v}}}}^h \in \hat{{V}}_0^h
\end{aligned}
\label{momcon}
\end{equation}
\end{mdframed}

Above, on the $i^\text{th}$ facet $\Gamma_{e_i}$ of the $e^\text{th}$ element $\Omega_e$, $\lambda$ is an indicator function that takes on a value of unity if the facet is an inflow facet and a value of zero if it is an outflow facet, i.e.:
\begin{equation}
\begin{aligned}
\lambda =
\begin{cases}
1 \hspace{5mm} \text{if } \overline{{\bf{u}}}^h \cdot {\bf{n}} <0\\
0 \hspace{5mm} \text{if } \overline{{\bf{u}}}^h \cdot {\bf{n}} \geq 0
\end{cases}
\end{aligned}
\end{equation}
These definitions result from the fact that we have discretized the incompressible RANS equations using an upwind treatment of the advective component of the flux, as discussed in  \cite{rhebergen_hybridizable_2017}.

The penalty constant $C_{pen}$ arises from the fact that we have discretized the diffusive component of the flux using the symmetric interior penalty method \cite{arnold_penalty}.  As is typical with interior penalty methods, the constant $C_{pen}$ must be chosen sufficiently large as to ensure the resulting semi-discrete formulation is energy stable \cite{cockburn_note}.  To arrive at an intelligent choice for $C_{pen}$, one typically turns to trace inequalities.  Following Warburton and Hesthaven \cite{penalty1}, it can be shown that the following inequality holds in the two-dimensional setting:
\begin{equation}
|| \nabla^S {\bf{v}^h}||_{\partial \Omega_e}^2 \leq  (k+1)(k+2) \frac{\text{Perimeter}(\Omega_e)}{\text{Area} (\Omega_e)} || \nabla^S {\bf{v}^h}||_{\Omega_e}^2
\label{p2}
\end{equation}
for all $\Omega_e \in \mathcal{T}$ and ${\bf{v}}^h \in V^h$, and a similar result holds in the three-dimensional setting.  If we define the mesh size $h_e$ to be equal to $\text{Area}(\Omega_e)/\text{Perimeter}(\Omega_e)$, it suffices to select $C_{pen} \geq (k+1)(k+2)$.  We later demonstrate this choice yields an energy stable method.  In all of our later computations, we have selected $C_{pen} = (k+1)(k+2)$.  Note that the definition $h_e = \text{Area}(\Omega_e)/\text{Perimeter}(\Omega_e)$ is somewhat non-standard.  It is related not to the diameter of the element but rather to the diameter of the incircle of the element. However, we have found this definition to be critical in ensuring stability in the presence of boundary layer meshes containing high-aspect ratio elements.\\
\\
We now present a consistency result for our semi-discrete formulation.\\
\\
\textbf{Proposition 1} (Consistency)\\
\\
\textit{The semi-discrete HDG method presented in Eqs. (\ref{mass}-\ref{momcon}) is consistent provided the exact solution $(\bf{u},p)$ of the incompressible RANS equations is sufficiently smooth.  That is, Eqs. (\ref{mass}-\ref{momcon}) hold if we replace $\left( \overline{{\bf{u}}}^h, \hat{\overline{{\bf{u}}}}^h, \overline{p}^h,  \hat{\overline{p}}^h \right)$ with $\left( \overline{\bf{u}}, \overline{\bf{u}}|_{\tilde{\Gamma}}, \overline{p},  \overline{p}|_{\tilde{\Gamma}} \right)$.}\\
\\
\textit{Proof.} Note that for, a sufficiently smooth exact solution $({\bf{u}},p)$, Eq. (\ref{NS_S}) are satisfied in a point-wise manner.  Using this, we show that each of Eqs. (\ref{mass}-\ref{momcon}) hold if we replace $\left( \overline{{\bf{u}}}^h, \hat{\overline{{\bf{u}}}}^h, \overline{p}^h,  \hat{\overline{p}}^h \right)$ with $\left( \overline{\bf{u}}, \overline{\bf{u}}|_{\tilde{\Gamma}}, \overline{p},  \overline{p}|^{}_{\tilde{\Gamma}} \right)$.  We begin with the continuity equation, that is, Eq. (\ref{mass}).  This holds trivially since $\overline{{\bf{u}}}$ is divergence-free and hence:
\begin{equation}
\begin{aligned}
\sum_e \int_{\Omega_e}\frac{1}{\rho}  \nabla \cdot \overline{{\bf{u}}} q^h d\Omega = 0  \ \hspace{10.00mm} \forall q^h \in Q^h
\end{aligned}
\end{equation}
We next consider the continuity conservativity equation, Eq. (\ref{masscon}).  This holds since:
\begin{equation}
\begin{aligned}
\sum_e \sum_i \int_{\Gamma_{e_i}}\frac{1}{\rho} \overline{{\bf{u}}} \cdot {\bf{n}} \hat{q}^h d\Gamma - \int_{\partial \Omega}\frac{1}{\rho} \overline{{\bf{u}}}|_{\tilde{\Gamma}} \cdot {\bf{n}} \hat{q}^h d\Gamma = \int_{\tilde{\Gamma} \backslash \partial \Omega} \frac{1}{\rho} \llbracket \overline{{\bf{u}}} \rrbracket \hat{q}^h d\Gamma = 0 \ \hspace{10.00mm} \forall \hat{q}^h \in \hat{Q}^h
\end{aligned}
\end{equation}
as $\llbracket \overline{{\bf{u}}} \rrbracket = 0$ on each interior facet $F \in \mathcal{F}_{int}$. We now continue with the momentum equation, Eq. (\ref{mom}).  Through reverse integration by parts, we have:
\begin{equation}
\begin{aligned}
&&&\int_{\Omega} \frac{\partial \overline{{\bf{u}}}}{\partial t} \cdot {\bf{v}}^h d \Omega -\sum_e \int_{\Omega_e} (\overline{{\bf{u}}} \otimes \overline{{\bf{u}}}) : \nabla {\bf{v}}^h d \Omega - \sum_e \int_{\Omega_e} \frac{1}{\rho} \overline{p} \textbf{I} :  \nabla {\bf{v}}^h d \Omega + \sum_e \int_{\Omega_e}  2 (\nu +\nu_T) \nabla^S \overline{{\bf{u}}}  :  \nabla^S {{\bf{v}}}^hd \Omega \\
&&&+\sum_e \sum_i \int_{\Gamma_{e_i}} (\overline{{\bf{u}}} \otimes \overline{{\bf{u}}}) : ({\bf{v}}^h \otimes {\bf{n}}) d\Gamma
+\sum_e \sum_i \int_{\Gamma_{e_i}} (\overline{{\bf{u}}}|_{\tilde{\Gamma}} - \overline{{\bf{u}}}) \otimes \lambda \overline{{\bf{u}}} : ({\bf{v}}^h \otimes {\bf{n}}) d\Gamma
+\sum_e \sum_i \int_{\Gamma_{e_i}}\frac{1}{\rho} \overline{p}|_{\tilde{\Gamma}} \textbf{I} : ({\bf{v}}^h \otimes {\bf{n}}) d\Gamma \\
&&&-\sum_e \sum_i \int_{\Gamma_{e_i}}2 (\nu +\nu_T) \nabla^S \overline{{\bf{u}}}: ({\bf{v}}^h \otimes {\bf{n}}) d\Gamma
 -\sum_e \sum_i \int_{\Gamma_{e_i}}{\frac{2C_{pen}}{h_e} (\nu + \nu_T)} (\overline{{\bf{u}}}|_{\tilde{\Gamma}} - \overline{{\bf{u}}})   \otimes {\bf{n}}: ({\bf{v}}^h \otimes {\bf{n}}) d\Gamma \\
&&& + \sum_e \sum_i \int_{\Gamma_{e_i}}2 (\nu+\nu_T)[ (\overline{{\bf{u}}}|_{\tilde{\Gamma}} - \overline{{\bf{u}}})  \otimes {\bf{n}}] :  \nabla^S {{\bf{v}}}^h d\Gamma
 - \int_{\Omega} {\bf{f}} \cdot {\bf{v}}^h d\Omega \\
&&&= \sum_e \int_{\Omega_e} \left[ 
\frac{\partial \overline{{\bf{u}}}}{\partial t}
+\nabla \cdot (\overline{{\bf{u}}} \otimes \overline{{\bf{u}}}) 
+ \frac{1}{\rho} \nabla \overline{p}
- \nabla \cdot (2 (\nu +\nu_T) \nabla^S \overline{{\bf{u}}})
 - {\bf{f}}
 \right] 
 \cdot {\bf{v}}^h d \Omega  = \vec{0} \hspace{10.00mm}
\forall {\bf{v}}^h \in {V}^h
\end{aligned}
\end{equation}
since:
\begin{equation}
\begin{aligned}
\frac{\partial \overline{{\bf{u}}}}{\partial t} +\nabla \cdot (\overline{{\bf{u}}} \otimes \overline{{\bf{u}}})  + \frac{1}{\rho} \nabla \overline{p} - \nabla \cdot (2 (\nu +\nu_T) \nabla^S \overline{{\bf{u}}}) - {\bf{f}}  = {\bf{0}} &&  \text{ in } \Omega \times (0,\infty)
\end{aligned}
\end{equation}
\\
We finish with the momentum conservativity equation, Eq. (\ref{momcon}).  This holds since:
\begin{equation}
\begin{aligned}
&&&-\sum_e \sum_i \int_{\Gamma_{e_i}} (\overline{{\bf{u}}} \otimes \overline{{\bf{u}}}) : (\hat{{{\bf{v}}}}^h \otimes {\bf{n}}) d\Gamma
-\sum_e \sum_i \int_{\Gamma_{e_i}} (\overline{{\bf{u}}}|_{\tilde{\Gamma}} - \overline{{\bf{u}}})  \otimes \lambda \overline{{\bf{u}}}^h : (\hat{{{\bf{v}}}}^h \otimes {\bf{n}}) d\Gamma \\
&&&-\sum_e \sum_i \int_{\Gamma_{e_i}}\frac{1}{\rho} \overline{p}|_{\tilde{\Gamma}} \textbf{I} : (\hat{{{\bf{v}}}}^h \otimes {\bf{n}}) d\Gamma
+\sum_e \sum_i \int_{\Gamma_{e_i}}2 (\nu +\nu_T) \nabla^S\overline{{\bf{u}}}:  (\hat{{{\bf{v}}}}^h \otimes {\bf{n}}) d\Gamma \\
&&&+\sum_e \sum_i \int_{\Gamma_{e_i}}{\frac{2C_{pen}}{h_e} (\nu + \nu_T)} (\overline{{\bf{u}}}|_{\tilde{\Gamma}} - \overline{{\bf{u}}})   \otimes {\bf{n}}:  (\hat{{{\bf{v}}}}^h \otimes {\bf{n}}) d\Gamma \\
&&&+\int_{\Gamma_N} (1-\lambda) (\overline{{\bf{u}}}|_{\tilde{\Gamma}} \cdot {\bf{n}}) (\overline{{\bf{u}}}|_{\tilde{\Gamma}} \cdot \hat{{{\bf{v}}}}^h) d \Gamma
+ \int_{\Gamma_N} {\bf{h}} \cdot \hat{{{\bf{v}}}}^h d \Gamma \\
&&&=\int_{\tilde{\Gamma}/\partial \Omega} \left( - \left( \overline{\bf{u}} \cdot \hat{\bf{v}}^h \right) \llbracket \overline{\bf{u}} \rrbracket + \llbracket 2(\nu + \nu_T) \nabla^s \bar{\bf{u}} \rrbracket \cdot \hat{\bf{v}}^h \right) d\Gamma \\
&&&- \int_{\Gamma_N} \left( \left[-\frac{1}{\rho}\overline{p} \textbf{I} +  2 (\nu+\nu_T) \nabla^S \overline{\bf{u}}\right] \cdot {\bf{n}}  - \text{min}(\overline{\bf{u}} \cdot {\bf{n}},0) \overline{\bf{u}} - {\bf{h}} \right) \cdot \hat{{\bf{v}}}^h d\Gamma = 0
\hspace{19mm }\forall \hat{{{\bf{v}}}}^h \in \hat{{V}}^h
\end{aligned}
\end{equation}
as $\llbracket \overline{{\bf{u}}} \rrbracket = 0$ and $\llbracket 2 (\nu + \nu_T) \nabla^S \overline{{\bf{u}}} \rrbracket = {\bf{0}}$ on each interior facet $F \in \mathcal{F}_{int}$ and:
\begin{equation}
\begin{aligned}
\left[-\frac{1}{\rho}\overline{p}\textbf{I} +  2 (\nu+\nu_T) \nabla^S \overline{\bf{u}}\right] \cdot {\bf{n}}  - \text{min}(\overline{\bf{u}} \cdot {\bf{n}},0) \overline{\bf{u}} - {\bf{h}} & = {\bf{0}} && \text{ on } \Gamma_N \times (0,\infty)
\end{aligned}
\end{equation}
This completes the proof.$\qed$\\

We next present a result illustrating our semi-discrete formulation conserves mass in a pointwise manner.

\textbf{Proposition 2} (Pointwise Mass Conservation)\\
\\
\textit{If  $\overline{{\bf{u}}}^h \in V^h$ and $\hat{\overline{{\bf{u}}}}^h \in \hat{V}^h$ satisfy Eqs. (\ref{mass}) and (\ref{masscon}), with $V_h$ and $\hat{V}^h$ defined in Eq .(\ref{spaces}), then:}\\
\begin{equation}
\nabla \cdot \overline{{\bf{u}}}^h = 0 \hspace{5mm} \forall {\bf{x}} \in \Omega_e,\hspace{5mm} \forall \Omega_e \in \mathcal{T}
\label{pointwise}
\end{equation}
\textit{and:}
\begin{equation}
\begin{aligned}
\llbracket \overline{{\bf{u}}}^h \rrbracket = 0 \hspace{5mm} &\forall {\bf{x}} \in F,\hspace{5mm} \forall F \in \mathcal{F}_{\text{int}}\\
\overline{{\bf{u}}}^h \cdot {\bf{n}} = \hat{\overline{{\bf{u}}}}^h \cdot {\bf{n}} \hspace{5mm} &\forall {\bf{x}} \in F,\hspace{5mm} \forall F \in \mathcal{F}_{\text{bdy}}
\end{aligned}
\end{equation}
\\
\textit{Proof.} The proof follows in the same manner as the proof of Proposition 1 in \cite{rhebergen_hybridizable_2017}, though we review the proof as pointwise mass conservation is a key feature of our HDG method.  From Eq. (\ref{mass}) it follows that:
\begin{equation}
0 = \int_{\Omega_e} q^h \nabla \cdot \overline{{\bf{u}}}^h d\Omega \hspace{5mm} \forall q^h \in P_{k-1}(\Omega_e),\hspace{5mm} \forall \Omega_e \in \mathcal{T}
\end{equation}
\\
Since $\nabla \cdot \overline{{\bf{u}}}^h \in P_{k-1}(\Omega_e)$ for every $\Omega_e \in \mathcal{T}$, we can take $q^h = \nabla \cdot \overline{{\bf{u}}}^h$ in the above equation, yielding $\int_{\Omega_e} ( \nabla \cdot \overline{{\bf{u}}}^h )^2 d\Omega = 0$ for every $\Omega_e \in \mathcal{T}$.  Thus $\nabla \cdot \overline{{\bf{u}}}^h \equiv 0$ in $\Omega$.\\
\\
From Eq. (\ref{masscon}) it follows that:
\begin{equation}
0 = \int_F  \llbracket \overline{{\bf{u}}}^h \rrbracket \hat{q}^h dS \hspace{5mm} \forall \hat{q}^h \in P_{k}(F),\hspace{5mm} \forall F \in \mathcal{F}_{\text{int}}
\end{equation}
Since $\llbracket \overline{{\bf{u}}}^h \rrbracket \in P_k(F)$ for all $F \in \mathcal{F}_{\text{int}}$, we can choose $\hat{q}^h = \llbracket \overline{{\bf{u}}}^h \rrbracket$, yielding $\int_{F} \llbracket \overline{{\bf{u}}}^h \rrbracket^2 d\Gamma = 0$ for all $F \in \mathcal{F}_{\text{int}}$.  Thus $\llbracket \overline{{\bf{u}}}^h \rrbracket \equiv 0$ on $\tilde{\Gamma}/\partial \Omega$.\\
\\
From Eq. (\ref{masscon}) it also follows that:
\begin{equation}
0 = \int_F  \left( \overline{{\bf{u}}}^h - \hat{\overline{{\bf{u}}}}^h \right) \cdot {\bf{n}} \hat{q}^h dS \hspace{5mm} \forall \hat{q}^h \in P_{k}(F),\hspace{5mm} \forall F \in \mathcal{F}_{\text{bdy}}
\end{equation}
Since $\left( \overline{{\bf{u}}}^h - \hat{\overline{{\bf{u}}}}^h \right) \cdot {\bf{n}} \in P_k(F)$ for all $F \in \mathcal{F}_{\text{bdy}}$, we can choose $\hat{q}^h = \left( \overline{{\bf{u}}}^h - \hat{\overline{{\bf{u}}}}^h \right) \cdot {\bf{n}}$, yielding $\int_{F} \left( \left( \overline{{\bf{u}}}^h - \hat{\overline{{\bf{u}}}}^h \right) \cdot {\bf{n}} \right)^2 d\Gamma = 0$ for all $F \in \mathcal{F}_{\text{bdy}}$.  Thus $\left( \overline{{\bf{u}}}^h - \hat{\overline{{\bf{u}}}}^h \right) \cdot {\bf{n}} \equiv 0$ on $\partial \Omega$.
$\qed$\\
\\
Proposition 2 is a stronger statement of mass conservation than general DG or HDG finite element methods, in which mass conservation is normally satisfied locally (element-wise) in an integral sense only.  There are two fundamental reasons why our semi-discrete formulation conserves mass in a pointwise manner.  First of all, the element-wise divergence operator maps the velocity space $V^h$ into the pressure space $Q^h$.  Second of all, the element-wise trace operator maps the velocity space $V^h$ into the pressure trace space $\hat{Q}^h$.  Any choice of velocity, pressure, velocity trace, and velocity pressure spaces that satisfy these two properties will result in a divergence-conforming HDG method.  It can easily be seen that the choice of approximation spaces presented here corresponds to a hybridization of a divergence-conforming DG method employing a Brezzi-Douglas-Marini velocity-pressure pair \cite{two_families}.\\
\\
The next result gives a global momentum balance law for our semi-discrete method.

\textbf{Proposition 3} (Global Momentum Balance)\\
\\
\textit{Let $(\overline{{\bf{u}}}^h, \hat{\overline{{\bf{u}}}}^h, \overline{p}^h, \hat{\overline{p}}^h) \in V^h \times \hat{V}^h \times Q^h \times \hat{Q}^h$ satisfy Eq.(\ref{mass} -\ref{momcon}).} 
\textit{If $\Gamma_D = \emptyset$:}
\begin{equation}
\begin{aligned}
\frac{d}{dt} &\int_{\Omega} \overline{{\bf{u}}}^h d\Omega =
\int_{\Omega} {\bf{f}} d\Omega -
\int_{\Gamma_N} (1-\lambda) (\hat{\overline{{\bf{u}}}}^h \cdot {\bf{n}}) \hat{\overline{{\bf{u}}}}^h  d \Gamma
- \int_{\Gamma_N} {\bf{h}}  d \Gamma
\end{aligned}
\label{momproof2}
\end{equation}
\textit{Proof.} The proof follows in the same manner as the proof of Proposition 2 in \cite{rhebergen_hybridizable_2017}. \qed\\
\\
The next result states that our semi-discrete method is globally energy stable.\\
\\
\textbf{Proposition 4} (Global Energy Stability)\\
\\
\textit{If $(\overline{{\bf{u}}}^h, \hat{\overline{{\bf{u}}}}^h, \overline{p}^h, \hat{\overline{p}}^h) \in V^h \times \hat{V}^h \times Q^h \times \hat{Q}^h$ satisfy Eqs.~(\ref{mass} -~\ref{momcon}), ${\bf{f}} = {\bf{0}}$, ${\bf{g}} = {\bf{0}}$, ${\bf{h}} = {\bf{0}}$, and $\nu_T \geq 0$:}
\begin{equation}
\frac{d}{dt}\sum_e \int_{\Omega_e} |\overline{{\bf{u}}}^h|^2 d\Omega_e \leq 0
\label{energy0}
\end{equation}
\textit{Proof}.  The proof follows in the same manner as the proof of Proposition 3 in \cite{rhebergen_hybridizable_2017}, though we review the proof to illustrate the role of the penalty constant $C_{pen}$. 
Setting ${\bf{v}}^h = \overline{{\bf{u}}}^h$ and $\hat{{\bf{v}}}^h = \hat{\overline{{\bf{u}}}}^h$ in Eqs. (\ref{mass}-\ref{momcon}), summing the two expressions together, recognizing that $\nabla \cdot \overline{\bf{u}}^h = 0$ in $\Omega$ and $\overline{{\bf{u}}}^h \cdot {\bf{n}} = \hat{\overline{{\bf{u}}}}^h  \cdot {\bf{n}}$ for all $F \in \mathcal{F}$, and exploiting that $\lambda \overline{{\bf{u}}}^h \cdot {\bf{n}} = (\overline{{\bf{u}}}^h \cdot {\bf{n}} - |\overline{{\bf{u}}}^h \cdot {\bf{n}}|)/2$, we obtain:

\begin{equation}
\begin{aligned}
\sum_e\frac{1}{2}&\int_{\Omega_e} \frac{\partial | \overline{{\bf{u}}}^h|^2}{\partial t} d\Omega
+\sum_e \sum_i \frac{1}{2} \int_{\Gamma_{e_i}} (\overline{{\bf{u}}}^h \cdot {\bf{n}})|\overline{{\bf{u}}}^h|^2 d\Gamma\\
-\sum_e \sum_i \frac{1}{2} &\int_{\Gamma_{e_i}}  (\overline{{\bf{u}}}^h \cdot {\bf{n}})|\hat{\overline{{\bf{u}}}}^h|^2d\Gamma
+\sum_e \sum_i \frac{1}{2} \int_{\Gamma_{e_i}} |\overline{{\bf{u}}}^h \cdot {\bf{n}} | | \overline{{\bf{u}}}^h - \hat{\overline{{\bf{u}}}}^h|^2 d\Gamma\\
+\sum_e&\int_{\Omega_e} 2(\nu+\nu_T)| \nabla^S \overline{{\bf{u}}}^h|^2d\Omega
+\sum_e \sum_i \int_{\Gamma_{e_i}} \frac{2C_{pen}}{h_e}(\nu+\nu_T)|\hat{\overline{{\bf{u}}}}^h - \overline{{\bf{u}}}^h|^2 d\Gamma\\
+4\sum_e \sum_i  &\int_{\Gamma_{e_i}} (\nu +\nu_T) (\nabla^S \overline{{\bf{u}}}^h \cdot {\bf{n}}) \cdot (\hat{\overline{{\bf{u}}}}^h - \overline{{\bf{u}}}^h)d\Gamma
+\int_{\Gamma_N} (1-\lambda) (\hat{\overline{{\bf{u}}}}^h \cdot {\bf{n}}) | \hat{\overline{{\bf{u}}}}^h|^2 d\Gamma\\
-\sum_e &\int_{\Omega_e} (\overline{{\bf{u}}}^h \otimes \overline{{\bf{u}}}^h) : \nabla \overline{{\bf{u}}}^h d\Omega = 0
\end{aligned}
\label{energy1}
\end{equation}
Since $\hat{\overline{{\bf{u}}}}^h$ is continuous across facets, $\llbracket \overline{{\bf{u}}}^h \rrbracket = 0$ on interior facets, and $\hat{\overline{{\bf{u}}}}^h \cdot {\bf{n}} = \overline{{\bf{u}}}^h \cdot {\bf{n}}$ on the domain boundary, the third integral on the left hand side of Eq.(\ref{energy1}) can be simplified as follows:
\begin{equation}
-\sum_e \sum_i \frac{1}{2} \int_{\Gamma_{e_i}}  (\overline{{\bf{u}}}^h \cdot {\bf{n}})|\hat{\overline{{\bf{u}}}}^h|^2d\Gamma =
-\frac{1}{2} \int_{\Gamma_N}(\hat{\overline{{\bf{u}}}}^h \cdot {\bf{n}}) |\hat{\overline{{\bf{u}}}}^h|^2 d\Gamma
\end{equation}
By the product rule, it holds that $-\overline{{\bf{u}}}^h \otimes \overline{{\bf{u}}}^h : \nabla \overline{{\bf{u}}}^h = -\nabla \cdot ((\overline{{\bf{u}}}^h \otimes \overline{{\bf{u}}}^h)\cdot \overline{{\bf{u}}}^h)/2$ since $\nabla \cdot \overline{{\bf{u}}}^h = 0$ on each element $\Omega_e$.  By the divergence theorem, it then holds that the last term on the left hand side of Eq. (\ref{energy1}) is equal to:
\begin{equation}
-\sum_e \int_{\Omega_e} (\overline{{\bf{u}}}^h \otimes \overline{{\bf{u}}}^h) : \nabla \overline{{\bf{u}}}^h d\Omega =
-\frac{1}{2}\sum_e \sum_i \int_{\Gamma_{e_i}} (\overline{{\bf{u}}}^h \cdot  {\bf{n}}) |\overline{{\bf{u}}}^h|^2 d\Gamma
\label{energy2}
\end{equation}

Combining Eqs. (\ref{energy1}-\ref{energy2}) yields:
\begin{equation}
\begin{aligned}
\sum_e\frac{1}{2}\int_{\Omega_e} \frac{\partial |\overline{{\bf{u}}}^h|^2}{\partial t} d\Omega =
- \sum_e \sum_i \frac{1}{2} \int_{\Gamma_{e_i}} |\overline{{\bf{u}}}^h \cdot {\bf{n}} | | \overline{{\bf{u}}}^h - \hat{\overline{{\bf{u}}}}^h|^2 d\Gamma \\
- \sum_e\int_{\Omega_e} 2(\nu+\nu_T)|\nabla^S \overline{{\bf{u}}}^h|^2d\Omega
- \sum_e \sum_i \int_{\Gamma_{e_i}} \frac{2C_{pen}}{h_e}(\nu+\nu_T)|\hat{\overline{{\bf{u}}}}^h - \overline{{\bf{u}}}^h|^2 d\Gamma\\
- 4\sum_e \sum_i  \int_{\Gamma_{e_i}} (\nu +\nu_T) (\nabla^S \overline{{\bf{u}}}^h \cdot {\bf{n}}) \cdot (\hat{\overline{{\bf{u}}}}^h - \overline{{\bf{u}}}^h)d\Gamma
- \frac{1}{2}\sum_e \sum_i \int_{\Gamma_{e_i}} (\hat{\overline{{\bf{u}}}}^h \cdot  {\bf{n}}) |\hat{\overline{{\bf{u}}}}^h|^2 d\Gamma
\end{aligned}
\label{energy3}
\end{equation}
where we have used that:
\begin{equation}
\int_{\Gamma_N} (1-\lambda) (\hat{\overline{{\bf{u}}}}^h \cdot {\bf{n}}) | \hat{\overline{{\bf{u}}}}^h|^2 d\Gamma
-\frac{1}{2} \int_{\Gamma_N}(\hat{\overline{{\bf{u}}}}^h \cdot {\bf{n}}) |\hat{\overline{{\bf{u}}}}^h|^2 d\Gamma =
\frac{1}{2}\sum_e \sum_i \int_{\Gamma_{e_i}} (\hat{\overline{{\bf{u}}}}^h \cdot  {\bf{n}}) |\hat{\overline{{\bf{u}}}}^h|^2 d\Gamma
\end{equation}
By the Cauchy-Schwarz inequality and Young's inequality, it holds that:
\begin{equation}
\begin{aligned}
& -4\sum_e \sum_i  \int_{\Gamma_{e_i}} (\nu +\nu_T) (\nabla^S \overline{{\bf{u}}}^h \cdot {\bf{n}}) \cdot (\hat{\overline{{\bf{u}}}}^h - \overline{{\bf{u}}}^h)d\Gamma \leq \\ & \sum_e  \int_{\partial \Omega_e} \frac{2h_e}{C_{pen}} (\nu+\nu_T)|\nabla^S \overline{{\bf{u}}}^h|^2d\Gamma + \sum_e \sum_i \int_{\Gamma_{e_i}} \frac{2C_{pen}}{h_e}(\nu+\nu_T)|\hat{\overline{{\bf{u}}}}^h - \overline{{\bf{u}}}^h|^2 d\Gamma
\label{energy4}
\end{aligned}
\end{equation}
and by the trace inequality, it holds that
\begin{equation}
\int_{\partial \Omega_e} (\nu+\nu_T)|\nabla^S \overline{{\bf{u}}}^h|^2d\Gamma \leq \frac{(k+1)(k+2)}{h_e} \int_{\Omega_e} (\nu+\nu_T)|\nabla^S \overline{{\bf{u}}}^h|^2d\Omega
\label{energy5}
\end{equation}
for each element $\Omega_e$.
Combining Eqs. (\ref{energy3}-\ref{energy5}) and exploiting the fact that $C_{pen} \geq (k+1)(k+2)$ yields:
\begin{equation}
\begin{aligned}
\sum_e\frac{1}{2}\int_{\Omega_e} \frac{\partial |\overline{{\bf{u}}}^h|^2}{\partial t} d\Omega \leq
- \sum_e \sum_i \frac{1}{2} \int_{\Gamma_{e_i}} |\overline{{\bf{u}}}^h \cdot {\bf{n}} | | \overline{{\bf{u}}}^h - \hat{\overline{{\bf{u}}}}^h|^2 d\Gamma 
- \frac{1}{2}\sum_e \sum_i \int_{\Gamma_{e_i}} (\hat{\overline{{\bf{u}}}}^h \cdot  {\bf{n}}) |\hat{\overline{{\bf{u}}}}^h|^2 d\Gamma \leq 0
\end{aligned}
\label{energy6}
\end{equation}
This completes the proof.
$\qed$

\section{Strong Form of the Advection-Diffusion Equation}
The HDG method that we have presented for the incompressible RANS equations applies to all linear eddy viscosity models.  With most eddy viscosity models, we must solve additional transport equations for so-called \textit{turbulence variables}.  For example, with a $k$-$\varepsilon$ model \cite{k_eps}, we solve transport equations for the turbulent kinetic energy $k$ and the turbulent dissipation $\varepsilon$, while with a Spalart-Allmaras model, we solve a single transport equation for the eddy visosity $\nu_T$ itself.  The transport equations for the turbulence variables typically take the form of \textit{advection-diffusion equations}.  As such, we present here HDG methods for a general advection-diffusion equation.

Once again consider a Lipschitz and bounded domain $\Omega \subset \mathbb{R}^d$, and let us partition the boundary of the domain $\partial \Omega$ into a Dirichlet boundary $\Gamma_D$ and a Neumann boundary $\Gamma_N$ such that $\partial \Omega = \overline{\Gamma_D \cup \Gamma_N}$ and $\Gamma_D \cap \Gamma_N = \emptyset$.  Let $f: \Omega \times (0,\infty) \rightarrow \mathbb{R}$, $g: \Gamma_D \times (0,\infty) \rightarrow \mathbb{R}$, $h: \Gamma_N \times (0,\infty) \rightarrow \mathbb{R}$, and $\phi_0: \Omega \rightarrow \mathbb{R}$.  Moreover, let $\kappa: \Omega \times (0,\infty) \rightarrow \mathbb{R}^+$, let $\overline{\bf{u}} : \Omega \times (0,\infty) \rightarrow \mathbb{R}^d$, and assume $\nabla \cdot \overline{\bf{u}} \equiv 0$.  The strong form of the advection-diffusion equation then reads:
\begin{mdframed}
\textbf{Strong Form of the Advection-Diffusion Equation}\\
\\
Find $\phi : \bar{\Omega} \times [0,\infty) \rightarrow \mathbb{R}$ such that:
\begin{equation}
\begin{aligned}
\frac{\partial \phi}{\partial t} + \overline{{\bf{u}}} \cdot \nabla \phi - \nabla \cdot (\kappa \nabla \phi) & = f  && \text{ in } \Omega \times (0,\infty) \label{adreq} \\
\phi & = g  && \text{ on } \Gamma_D \times (0,\infty)\\
\overline{{\bf{u}}} \phi \cdot {\bf{n}} - \kappa \nabla \phi  \cdot {\bf{n}}  - \text{max}(\overline{{\bf{u}}} \cdot {\bf{n}},0)\phi &= h  && \text{ on } \Gamma_N \times (0,\infty)\\
\phi(\cdot,0) & = \phi_0  && \text{ in } \Omega\\
\end{aligned}
\end{equation}
\end{mdframed}
Note that on inflow parts of  $\Gamma_N$ ($\overline{{\bf{u}}} \cdot {\bf{n}} < 0$ ) we impose the total flux, i.e., $\overline{{\bf{u}}} \phi \cdot {\bf{n}} - \kappa \nabla \phi  \cdot {\bf{n}} =h$, and on outflow parts of $\Gamma_N$ ( $\overline{{\bf{u}}} \cdot {\bf{n}}  \geq 0$), only the diffusive part of the flux is prescribed, i.e., $- \kappa \nabla \phi  \cdot {\bf{n}}=h$.\\
\\
From the above strong form, we see that $\phi$ denotes the transported scalar, $f$ denotes the forcing, $g$ denotes the applied Dirichlet boundary condition, $h$ denotes the applied Neumann boundary condition, $\phi_0$ denotes the initial condition, $\kappa$ denotes the diffusivity, and $\overline{\bf{u}}$ denotes the advection velocity.  We have utilized the notation $\overline{*}$ for the advection velocity as turbulence variables are advected with the mean velocity field.  At this stage, however, we consider the advection velocity to be known and fixed.  Later, we will couple the incompressible RANS equations with a particular eddy viscosity model, namely the Spalart-Allmaras model, in which case the advection velocity is also an unknown.

\section{Semi-Discrete HDG Formulation for the Advection-Diffusion Equation}
\label{sec:semidiscrete_AD}
We now construct an HDG method for the advection-diffusion equation.  As we did with the incompressible RANS equations, we first discretize in space, and later, we discretize in time.  In our presentation, we also employ the same notation introduced earlier for the incompressible RANS equations, except that we consider the new finite element spaces:
\begin{equation}
\begin{aligned}
&W^h:= \left\{ w^h \in L^2(\Omega): \left.w^h\right|_{\Omega_e} \in  P_k(\Omega_e) \hspace{2mm} \forall   \Omega_e \in \mathcal{T} \right\} \\
&\hat{W}^h:= \left\{ \hat{w}^h \in L^2(\tilde\Gamma): \left.\hat{w}^h\right|_{F} \in  P_k(F)  \hspace{2mm} \forall   F \in \mathcal{F} \right\} \\
&\hat{W}_g^h:= \left\{ \hat{w}^h \in \hat{W}^h: \hat{w}^h = g \text{ on } \Gamma_D
\right\} \\
&\hat{W}_0^h:= \left\{ \hat{w}^h \in \hat{W}^h: \hat{w}^h = 0 \text{ on } \Gamma_D
\right\} \\
\end{aligned}
\end{equation}
We will approximate the transported scalar over element interiors using the space $W^h$ and over the mesh skeleton using the space $\hat{W}^h$.  We henceforth refer to $W^h$ as the (discrete) scalar space and $\hat{W}^h$ are the (discrete) scalar trace space.  With these spaces in hand, our semi-discrete HDG formulation for the advection-diffusion equation takes the following form, wherein advective fluxes are treated using upwinding, diffusive terms are treated using the symmetric interior penalty method, and an artificial diffusivity is introduced to prevent spurious oscillations in the presence of sharp layers.
\\

\begin{mdframed}
\textbf{Semi-Discrete Formulation of the Advection-Diffusion Equation}\\
\\
Find $\left( \phi^h, \hat{\phi}^h \right) \in W^h \times \hat{W}_g^h$ such that:\\
\\
\noindent \textit{Advection-Diffusion Equation}
\begin{equation}
\begin{aligned}
&\int_{\Omega} \frac{\partial \phi}{\partial t}^h w^h d\Omega
- \sum_e \int_{\Omega_e}\overline{{\bf{u}}} \phi^h \cdot \nabla w^h d\Omega
+\sum_e \int_{\Omega_e} (\kappa+\kappa_{\text{DC}}) \nabla \phi^h \cdot \nabla w^h d \Omega \\
 \sum_e \sum_i &\int_{\Gamma_{e_i}} \overline{{\bf{u}}} \phi^h \cdot w^h {\bf{n}} d\Gamma
 +\sum_e \sum_i \int_{\Gamma_{e_i}} (\hat{\phi}^h - \phi^h) \lambda \overline{{\bf{u}}} \cdot w^h {\bf{n}} d \Gamma \\
 - \sum_e \sum_i &\int_{\Gamma_{e_i}} \kappa \nabla \phi^h \cdot w^h {\bf{n}} d\Gamma
 + \sum_e \sum_i \int_{\Gamma_{e_i}} \kappa (\hat{\phi}^h - \phi^h) {\bf{n}} \cdot \nabla w^h d \Gamma\\
 - \sum_e \sum_i &\int_{\Gamma_{e_i}} \frac{C_{\text{pen}}}{h_e} \kappa (\hat{\phi}^h - \phi^h){\bf{n}} \cdot w^h {\bf{n}} d \Gamma
 -\int_{\Omega} f w^h d \Omega = 0
 \hspace{30mm} \forall w^h \in W^h
\end{aligned}
\label{sclrone}
\end{equation}
\newpage
\noindent \textit{Advection-Diffusion Conservativity Condition}
\begin{equation}
\begin{aligned}
\sum_e \sum_i &\int_{\Gamma_{e_i}} \overline{{\bf{u}}} \phi^h \cdot \hat{w}^h {\bf{n}} d \Gamma
+ \sum_e \sum_i \int_{\Gamma_{e_i}} (\hat{\phi}^h - \phi^h) \lambda \overline{{\bf{u}}} \cdot \hat{w}^h {\bf{n}} d \Gamma \\
-\sum_e \sum_i &\int_{\Gamma_{e_i}} \kappa  \nabla \phi^h \cdot \hat{w}^h {\bf{n}} d \Gamma
-\sum_e \sum_i \int_{\Gamma_{e_i}} \frac{C_{\text{pen}}}{h_e} \kappa (\hat{\phi}^h - \phi^h ) {\bf{n}} \cdot \hat{w}^h {\bf{n}} d \Gamma \\
-\sum_e \sum_i &\int_{\Gamma_{e_i}} (1- \lambda) \hat{\phi}^h {\bf{n}} \cdot \overline{{\bf{u}}} \hat{w}^h d \Gamma - \int_{\Gamma_N} h \hat{w}^h d \Gamma = 0
 \hspace{38.5mm} \forall \hat{w}^h \in \hat{W}_0^h
\end{aligned}
\label{sclrtwo}
\end{equation}
\end{mdframed}

The penalty parameter $C_{pen}$ appearing in our semi-discrete formulation arises from the fact that we have discretized the diffusive component of the flux using the symmetric interior penalty method.  As was the case for the incompressible RANS equations, the penalty parameter must be chosen sufficiently large to ensure the resulting semi-discrete formulation is energy stable, and it suffices to choose $C_{pen} \geq (k+1)(k+2)$.

We have also included an artificial diffusivity parameter $\kappa_{\text{DC}}$ in our semi-discrete formulation to introduce numerical dissipation in regions of the domain exhibiting sharp gradients.  Provided the artificial diffusivity parameter is chosen sufficiently large, this ensures our numerical approximation of the transported scalar does not experience spurious oscillations.  This is especially important when transporting turbulence variables which must remain non-negative in the domain.  In what follows, we assume the artificial diffusivity parameter $\kappa_{\text{DC}}$ depends on the residual of the advection-diffusion equation.  In our later numerical experiments, we employ the artificial diffusivity parameter associated with the $YZ\beta$ method \cite{yzbeta2,yzbeta}:
\begin{equation}
\begin{aligned}
\kappa_{\text{DC}} = \left| Y^{-1} Z \right| \left( \left| Y^{-1} \nabla \phi^h \right|^2 \right)^{\beta/2 -1} \left( \frac{h_{\text{DC}}}{2}\right)^{\beta}
\label{yzb4}
\end{aligned}
\end{equation}
where $Y$ is a reference value of the scalar field $\phi^h$, $Z = \frac{\partial \phi^h}{\partial t} + \overline{{\bf{u}}} \cdot \nabla \phi^h - \nabla \cdot (\kappa \nabla \phi^h) - f$ is the point-wise residual of the advection-diffusion equation, $h_{\text{DC}}$ is a local element length scale which is precisely defined in \cite{yzbeta2}, and $\beta$ is a parameter (which we later choose to be equal to two) that influences the smoothness of layers.\\
\\
We now present a consistency result for our semi-discrete formulation.\\
\\
\textbf{Proposition 4} (Consistency)\\
\\
\textit{The semi-discrete HDG method presented in Eqs. (\ref{sclrone}-\ref{sclrtwo}) is consistent provided the exact solution $\phi$ of the advection-diffusion equation is sufficiently smooth.  That is, Eqs. (\ref{sclrone}-\ref{sclrtwo}) hold if we replace $\left( \phi^h, \hat{\phi}^h \right)$  with $\left( \phi, \left.\phi \right|_{\tilde{\Gamma}}\right)$}.\\
\\
\textit{Proof.}  Note that for, a sufficiently smooth exact solution $\phi$, Eq. (\ref{adreq}) is satisfied in a point-wise manner.  Using this, we show that each of Eqs. (\ref{sclrone}-\ref{sclrtwo}) hold if we replace $\left( \phi^h, \hat{\phi}^h\right)$  with $\left( \phi, \left.\phi \right|_{\tilde{\Gamma}} \right)$.  We begin with the advection-diffusion equation, that is,  Eq. (\ref{sclrone}).  Through reverse integration by parts we have: 

\begin{equation} for 
\begin{aligned}
&\int_{\Omega} \frac{\partial \phi}{\partial t} w^h d\Omega
- \sum_e \int_{\Omega_e}\overline{{\bf{u}}} \phi \cdot \nabla w^h d\Omega
+\sum_e \int_{\Omega_e} (\kappa+\kappa_{\text{DC}}) \nabla \phi \cdot \nabla w^h d \Omega \\
 \sum_e \sum_i &\int_{\Gamma_{e_i}} \overline{{\bf{u}}} \phi \cdot w^h {\bf{n}} d\Gamma
 +\sum_e \sum_i \int_{\Gamma_{e_i}} (\left.\phi \right|_{\tilde{\Gamma}} - \phi) \lambda \overline{{\bf{u}}} \cdot w^h {\bf{n}} d \Gamma \\
 - \sum_e \sum_i &\int_{\Gamma_{e_i}} \kappa \nabla \phi^h \cdot w^h {\bf{n}} d\Gamma
 + \sum_e \sum_i \int_{\Gamma_{e_i}} \kappa (\left.\phi \right|_{\tilde{\Gamma}} - \phi) {\bf{n}} \cdot \nabla w^h d \Gamma\\
 - \sum_e \sum_i &\int_{\Gamma_{e_i}} \frac{C_{\text{pen}}}{h_e} \kappa (\left.\phi \right|_{\tilde{\Gamma}} - \phi){\bf{n}} \cdot w^h {\bf{n}} d \Gamma
 -\int_{\Omega} f w^h d \Omega \\
  =\sum_e &\int_{\Omega_e} \left[ \frac{\partial \phi}{\partial t} + \overline{{\bf{u}}} \cdot \nabla \phi - \nabla \cdot (\kappa \nabla \phi) - f \right] w^h d\Omega = 0
 \hspace{30mm} \forall w^h \in W^h
\end{aligned}
\end{equation}
\\
since: 
\begin{equation}
\frac{\partial \phi}{\partial t} + \overline{{\bf{u}}} \cdot \nabla \phi - \nabla \cdot (\kappa \nabla \phi) - f = 0  \hspace{10mm} \text{   in } \Omega \times (0,\infty)
\label{adrsf}
\end{equation}
\\
Note that above $\kappa_{DC}$ is identically zero as it is a function of the residual which is zero.  We finish with the advection-diffusion conservativity equation, Eq. (\ref{sclrtwo}).  This holds since:
\begin{equation}
\begin{aligned}
\sum_e \sum_i &\int_{\Gamma_{e_i}} \overline{{\bf{u}}} \phi \cdot \hat{w}^h {\bf{n}} d \Gamma
+ \sum_e \sum_i \int_{\Gamma_{e_i}} ( \left.\phi \right|_{\tilde{\Gamma}} - \phi) \lambda \overline{{\bf{u}}} \cdot \hat{w}^h {\bf{n}} d \Gamma \\
-\sum_e \sum_i &\int_{\Gamma_{e_i}} \kappa  \nabla \phi \cdot \hat{w}^h {\bf{n}} d \Gamma
-\sum_e \sum_i \int_{\Gamma_{e_i}} \frac{C_{\text{pen}}}{h_e} \kappa ( \left.\phi \right|_{\tilde{\Gamma}} - \phi ) {\bf{n}} \cdot \hat{w}^h {\bf{n}} d \Gamma \\
-\sum_e \sum_i &\int_{\Gamma_{e_i}} (1- \lambda)  \left.\phi \right|_{\tilde{\Gamma}} {\bf{n}} \cdot \hat{\overline{{\bf{u}}}} \hat{w}^h d \Gamma - \int_{\Gamma_N} h \hat{w}^h d \Gamma\\
=&\int_{\tilde{\Gamma}/\partial \Omega} \left( \overline{\bf{u}} \hat{{w}}^h  \llbracket {\bf{\phi}} \rrbracket + \llbracket \kappa \nabla \phi \rrbracket \hat{{w}}^h \right) d\Gamma\\ 
-&\int_{\Gamma_N} \left( \left[\overline{{\bf{u}}} \phi - \kappa \nabla \phi\right]\cdot {\bf{n}}  - \text{max}(\overline{{\bf{u}}} \cdot {\bf{n}},0)\phi - h \right) \hat{w}^hd\Gamma = 0
\hspace{26mm} \forall \hat{w}_h \in \hat{W}_h
\end{aligned}
\end{equation}
as  $\llbracket {\bf{\phi}} \rrbracket = 0$ and $\llbracket \kappa \nabla \phi \rrbracket = 0$ on each interior facet $F \in \mathcal{F}_{int}$ and:
\begin{equation}
\begin{aligned}
\left[\overline{{\bf{u}}} \phi - \kappa \nabla \phi\right]\cdot {\bf{n}}  - \text{max}(\overline{{\bf{u}}} \cdot {\bf{n}},0)\phi - h = 0 \hspace{10pt} \text{ on  } \Gamma_{N} \times (0,\infty)
\end{aligned}
\end{equation}

This completes the proof.$\qed$\\

\section{Semi-Discrete HDG Formulation for the Reynolds Averaged Navier-Stokes Equations with the Spalart-Allmaras Turbulence Model}
We have now constructed semi-discrete HDG formulations for (i) the incompressible RANS equations given a non-negative turbulent viscosity and (ii) the advection-diffusion equation.  In this section, we present a semi-discrete HDG formulation for the incompressible RANS equations and a particular linear eddy viscosity model, namely the Spalart-Allmaras one equation turbulence model.  With the Spalart-Allmaras model, one solves for a \textit{working viscosity} $\tilde{\nu}$ using the model transport equation:
\begin{equation}
\begin{aligned}
\frac{\partial \tilde{\nu} }{\partial t} + \overline{{\bf{u}}} \cdot \nabla \tilde{\nu} - \nabla \cdot \frac{1}{\sigma} \left[ (\nu + \tilde{\nu} ) \nabla \tilde{\nu} \right] - C_{b1} \tilde{S} \tilde{\nu} + C_{w1} f_{w} \frac{\tilde{\nu}^2}{d^2} - \frac{C_{b2}}{\sigma} \nabla \tilde{\nu} \cdot \nabla \tilde{\nu} = 0
\end{aligned}
\end{equation}
and then the eddy viscosity $\nu_T$ is obtained from the working viscosity $\tilde{\nu}$ through a relation of the form:
\begin{equation}
\begin{aligned}
\nu_T = \tilde{\nu} f_{v1}
\end{aligned}
\end{equation}
The precise form for the functions $\tilde{S}$, $d$, $f_w$, and $f_{v1}$ as well as typical values for the model constants $\sigma$, $C_{b1}$, $C_{b2}$, $C_{w1}$, and $C_{v1}$ may be found in \cite{spalart_one-equation_1992}.  Note that the model transport equation for the working viscosity $\tilde{\nu}$ is an advection-diffusion equation of the form:
\begin{equation}
\frac{\partial \phi}{\partial t} + \overline{{\bf{u}}} \cdot \nabla \phi - \nabla \cdot (\kappa \nabla \phi) = f
\end{equation}
where $\phi = \tilde{\nu}$, $\kappa = \frac{\nu + \tilde{\nu}}{\sigma}$, and $f = C_{b1} \tilde{S} \tilde{\nu} - C_{w1} f_{w} \frac{\nu^2}{d^2} + \frac{C_{b2}}{\sigma} \nabla \tilde{\nu} \cdot \nabla \tilde{\nu}$.  As such, we can simply replace $\phi$, $\kappa$, and $f$ in our previously presented semi-discrete HDG formulation for the advection-diffusion equation to arrive at a semi-discrete HDG formulation for the Spalart-Allmaras turbulence model.  With this in mind, our semi-discrete HDG formulation for the incompressible RANS equations with the Spalart-Allmaras turbulence model is as follows:


\begin{mdframed}
\textbf{Semi-Discrete HDG Formulation for the Incompressible RANS Equations}\\
\\
Find $\left( \overline{{\bf{u}}}^h, \hat{\overline{{\bf{u}}}}^h, \overline{p}^h,  \hat{\overline{p}}^h, \tilde{\nu}^h, \hat{\tilde{\nu}}^h \right) \in V^h \times \hat{V}^h_g \times Q^h \times \hat{Q}^h \times W^h \times \hat{W}^h$ such that:\\
\\
\noindent \textit{Continuity Equation}
\begin{equation}
\begin{aligned}
\sum_e \int_{\Omega_e}\frac{1}{\rho}  \nabla \cdot \overline{{\bf{u}}}^h q^h d\Omega = 0 \hspace{81mm} \forall q^h \in Q^h
\end{aligned}
\end{equation}
\noindent \textit{Continuity Conservativity Condition}
\begin{equation}
\begin{aligned}
\sum_e \sum_i \int_{\Gamma_{e_i}}\frac{1}{\rho} \overline{{\bf{u}}}^h \cdot {\bf{n}} \hat{q}^h d\Gamma - \int_{\partial \Omega}\frac{1}{\rho} \hat{\overline{{\bf{u}}}}^h \cdot {\bf{n}} \hat{q}^h d\Gamma = 0 \hspace{46mm} \forall \hat{q}^h \in \hat{Q}^h
\end{aligned}
\end{equation}
\noindent \textit{Momentum Equation}
\begin{equation}
\begin{aligned}
&\int_{\Omega} \frac{\partial \overline{{\bf{u}}}^h}{\partial t} \cdot {\bf{v}}^h d \Omega \\
-\sum_e &\int_{\Omega_e} (\overline{{\bf{u}}}^h \otimes \overline{{\bf{u}}}^h) : \nabla {\bf{v}}^h d \Omega
- \sum_e \int_{\Omega_e} \frac{1}{\rho} \overline{p}^h \textbf{I} :  \nabla {\bf{v}}^h d \Omega
+ \sum_e \int_{\Omega_e}  2 (\nu +\nu_T(\tilde{\nu}^h)) \nabla^S \overline{{\bf{u}}}^h  :  \nabla^S {{\bf{v}}}^hd \Omega \\
+\sum_e \sum_i &\int_{\Gamma_{e_i}} (\overline{{\bf{u}}}^h \otimes \overline{{\bf{u}}}^h) : ({\bf{v}}^h \otimes {\bf{n}}) d\Gamma
+\sum_e \sum_i \int_{\Gamma_{e_i}} (\hat{\overline{{\bf{u}}}}^h - \overline{{\bf{u}}}^h) \otimes \lambda \overline{{\bf{u}}}^h : ({\bf{v}}^h \otimes {\bf{n}}) d\Gamma\\
+\sum_e \sum_i &\int_{\Gamma_{e_i}}\frac{1}{\rho} \hat{\overline{p}}^h \textbf{I} : ({\bf{v}}^h \otimes {\bf{n}}) d\Gamma
-\sum_e \sum_i \int_{\Gamma_{e_i}}2 (\nu +\nu_T(\tilde{\nu}^h)) \nabla^S \overline{{\bf{u}}}^h: ({\bf{v}}^h \otimes {\bf{n}}) d\Gamma \\
 -\sum_e \sum_i &\int_{\Gamma_{e_i}}{\frac{2C_{pen}}{h_e} (\nu + \nu_T(\tilde{\nu}^h))} (\hat{\overline{{\bf{u}}}}^h - \overline{{\bf{u}}}^h)   \otimes {\bf{n}}: ({\bf{v}}^h \otimes {\bf{n}}) d\Gamma \\
 + \sum_e \sum_i &\int_{\Gamma_{e_i}}2 (\nu+\nu_T(\tilde{\nu}^h))[ (\hat{\overline{{\bf{u}}}}^h - \overline{{\bf{u}}}^h)  \otimes {\bf{n}}] :  \nabla^S {{\bf{v}}}^h d\Gamma
 - \int_{\Omega} {\bf{f}} \cdot {\bf{v}}^h d\Omega = 0 \hspace{20.5mm}
\forall {\bf{v}}^h \in {V}^h
\end{aligned}
\end{equation}
\noindent \textit{Momentum Conservativity Condtion}
\begin{equation}
\begin{aligned}
-\sum_e \sum_i &\int_{\Gamma_{e_i}} (\overline{{\bf{u}}}^h \otimes \overline{{\bf{u}}}^h) : (\hat{{{\bf{v}}}}^h \otimes {\bf{n}}) d\Gamma
-\sum_e \sum_i \int_{\Gamma_{e_i}} (\hat{\overline{{\bf{u}}}}^h - \overline{{\bf{u}}}^h)  \otimes \lambda \overline{{\bf{u}}}^h : (\hat{{{\bf{v}}}}^h \otimes {\bf{n}}) d\Gamma \\
-\sum_e \sum_i &\int_{\Gamma_{e_i}}\frac{1}{\rho} \hat{\overline{p}}^h \textbf{I} : (\hat{{{\bf{v}}}}^h \otimes {\bf{n}}) d\Gamma
+\sum_e \sum_i \int_{\Gamma_{e_i}}2 (\nu +\nu_T(\tilde{\nu}^h)) \nabla^S\overline{{\bf{u}}}^h:  (\hat{{{\bf{v}}}}^h \otimes {\bf{n}}) d\Gamma \\
+\sum_e \sum_i &\int_{\Gamma_{e_i}}{\frac{2C_{pen}}{h_e} (\nu + \nu_T(\tilde{\nu}^h))} (\hat{\overline{{\bf{u}}}}^h - \overline{{\bf{u}}}^h)   \otimes {\bf{n}}:  (\hat{{{\bf{v}}}}^h \otimes {\bf{n}}) d\Gamma \\
+&\int_{\Gamma_N} (1-\lambda) (\hat{\overline{{\bf{u}}}}^h \cdot {\bf{n}}) (\hat{\overline{{\bf{u}}}}^h \cdot \hat{{{\bf{v}}}}^h) d \Gamma
+ \int_{\Gamma_N} {\bf{h}} \cdot \hat{{{\bf{v}}}}^h d \Gamma = 0
\hspace{29mm }\forall \hat{{{\bf{v}}}}^h \in \hat{{V}}^h_0
\end{aligned}
\end{equation}
\noindent \textit{Turbulence Model Equation}
\begin{equation}
\begin{aligned}
&\int_{\Omega} \frac{\partial \tilde{\nu}^h}{\partial t} w^h d\Omega
- \sum_e \int_{\Omega_e}\overline{{\bf{u}}}^h \tilde{\nu}^h \cdot \nabla w^h d\Omega
+\sum_e \int_{\Omega_e} \left(\frac{1}{\sigma}\left(\nu + \tilde{\nu}^h\right)+\kappa_{\text{DC}}\right) \nabla \tilde{\nu}^h \cdot \nabla w^h d \Omega \\
 \sum_e \sum_i &\int_{\Gamma_{e_i}} \hat{\overline{{\bf{u}}}}^h \tilde{\nu}^h \cdot w^h {\bf{n}} d\Gamma
 +\sum_e \sum_i \int_{\Gamma_{e_i}} (\hat{\tilde{\nu}}^h - \tilde{\nu}^h) \lambda \overline{{\bf{u}}} \cdot w^h {\bf{n}} d \Gamma \\
 - \sum_e \sum_i &\int_{\Gamma_{e_i}} \frac{1}{\sigma}\left(\nu + \tilde{\nu}^h\right) \nabla \tilde{\nu}^h \cdot w^h {\bf{n}} d\Gamma
 + \sum_e \sum_i \int_{\Gamma_{e_i}} \frac{1}{\sigma}\left(\nu + \tilde{\nu}^h\right) (\hat{\tilde{\nu}}^h - \tilde{\nu}^h) {\bf{n}} \cdot \nabla w^h d \Gamma\\
 - \sum_e \sum_i &\int_{\Gamma_{e_i}} \frac{C_{\text{pen}}}{h_e} \frac{1}{\sigma}\left(\nu + \tilde{\nu}^h\right) (\hat{\tilde{\nu}}^h - \tilde{\nu}^h){\bf{n}} \cdot w^h {\bf{n}} d \Gamma\\ 
 - &\int_{\Omega} \left( C_{b1} \tilde{S} \tilde{\nu}- C_{w1} f_{w} \frac{\nu^2}{d^2} + \frac{C_{b2}}{\sigma} \nabla \tilde{\nu} \cdot \nabla \tilde{\nu} \right) w^h d \Omega = 0
 \hspace{35.5mm} \forall w^h \in W^h
\end{aligned}
\end{equation}
\newpage
\noindent \textit{Turbulence Model Conservativity Condition}
\begin{equation}
\begin{aligned}
\sum_e \sum_i &\int_{\Gamma_{e_i}} \hat{\overline{{\bf{u}}}}^h \tilde{\nu}^h \cdot \hat{w}^h {\bf{n}} d \Gamma
+ \sum_e \sum_i \int_{\Gamma_{e_i}} (\hat{\tilde{\nu}}^h - \tilde{\nu}^h) \lambda \hat{\overline{{\bf{u}}}}^h \cdot \hat{w}^h {\bf{n}} d \Gamma \\
-\sum_e \sum_i &\int_{\Gamma_{e_i}} \frac{1}{\sigma}\left(\nu + \tilde{\nu}^h\right) \nabla \tilde{\nu}^h \cdot \hat{w}^h {\bf{n}} d \Gamma
-\sum_e \sum_i \int_{\Gamma_{e_i}} \frac{C_{\text{pen}}}{h_e} \frac{1}{\sigma}\left(\nu + \tilde{\nu}^h\right) (\hat{\tilde{\nu}}^h - \tilde{\nu}^h ) {\bf{n}} \cdot \hat{w}^h {\bf{n}} d \Gamma \\
-\sum_e \sum_i &\int_{\Gamma_{e_i}} (1- \lambda) \hat{\tilde{\nu}}^h {\bf{n}} \cdot \hat{\overline{{\bf{u}}}}^h \hat{w}^h d \Gamma - \int_{\Gamma_N} h \hat{w}^h d \Gamma = 0
 \hspace{40.25mm} \forall \hat{w}^h \in \hat{W}^h_0
\end{aligned}
\end{equation}
\end{mdframed}

By construction the above semi-discrete HDG formulation is consistent.  Moreover, the formulation yields a point-wise divergence-free velocity field, it conserves momentum globally, and, provided the formulation yields a non-negative eddy viscosity, it is energy stable.  These properties directly from the results obtained in Sections \ref{sec:semidiscrete_RANS} and \ref{sec:semidiscrete_AD}.


\section{Fully-Discrete HDG Formulation for the Reynolds Averaged Navier-Stokes Equations with the Spalart-Allmaras Turbulence Model}
Up to this point, we have discretized in space but not in time.  To discretize in time, we employ a staggered version of the so-called Generalized-$\alpha$ method.  This method was originally developed for structural mechanics by Chung and Hulbert in \cite{genalpha1} and then later extended to fluid mechanics in \cite{genalpha2} by Jansen et al.  To begin, let $\bf{U}$, $\bf{P}$, and $\dot{\bf{U}}$ denote the degree-of-freedom vectors associated with the interior velocity, pressure, and acceleration fields, the interior acceleration field, let $\hat{\bf{U}}$ and $\hat{\bf{P}}$ denote the degree-of-freedom vectors associated with the trace velocity and pressure fields, let $\bf{\Phi}$ and $\dot{\bf{\Phi}}$ denote the degree-of-freedom vectors associated with the interior working viscosity field and its time derivative, and let $\hat{\bf{\Phi}}$ denote the degree-of-freedom vector associated with the trace working viscosity field.  With the above notation in hand, we can write the semi-discrete HDG formulation for the incompressible RANS equations in terms of a system of differential-algebraic equations of the form:
\begin{align}
{\bf{R}}_{\text{flow}}(\dot{{\bf{U}}},{\bf{U}},{\bf{P}},\hat{{\bf{U}}},\hat{{\bf{P}}},{\bf{\Phi}}) &= {\bf{0}} \label{eq:residual_flow}\\
{\bf{R}}_{\text{scalar}}(\dot{{\bf{\Phi}}},{\bf{\Phi}},\hat{{\bf{\Phi}}},{\bf{U}}) &= {\bf{0}}  \label{eq:residual_scalar}
\end{align}
where Eq. \eqref{eq:residual_flow} contains the residual equations associated with the continuity equation, the continuity conservativity condition, the momentum equation, and the momentum conservativity condition while Eq. \eqref{eq:residual_scalar} contains the residual equations associated with the turbulence model equation and the turbulence model conservativity condition.  To solve the above system of equations at the $(n+1)^{\text{st}}$ time-step, we employ the Generalized-$\alpha$ method, freezing the working viscosity field to that associated with the $n^{\text{th}}$ step in Eq. \eqref{eq:residual_flow} and freezing the velocity field to that associated with the $n^{\text{th}}$ step in Eq. \eqref{eq:residual_scalar}.  This yields the following system of nonlinear algebraic equations:
\begin{align}
{\bf{R}}_{\text{flow}}(\dot{{\bf{U}}}_{n+\alpha_m},{\bf{U}}_{n+\alpha_f},{\bf{P}}_{n+1},\hat{{\bf{U}}}_{n+1},\hat{{\bf{P}}}_{n+1},{\bf{\Phi}}_{n}) &= {\bf{0}} \label{eq:genalpha_flow}\\
{\bf{R}}_{\text{scalar}}(\dot{{\bf{\Phi}}}_{n+\alpha_m},{\bf{\Phi}}_{n+\alpha_f},\hat{{\bf{\Phi}}}_{n+1},{\bf{U}}_n) &= {\bf{0}}  \label{eq:genalpha_scalar}
\end{align}
where $\alpha_m$ and $\alpha_f$ are determined from a free parameter $\rho_{\infty} \in [0,1]$ via:
\begin{equation}
\alpha_m = \frac{1}{2} \left( \frac{3 - \rho_{\infty}}{1+\rho_{\infty}} \right), \hspace{10mm} \alpha_f = \frac{1}{1+\rho_{\infty}}
\end{equation}
For a linear model problem, it has been shown the Generalized-$\alpha$ method annihilates the highest frequency in one time-step if $\rho_{\infty} = 0$ and preserves the highest frequency if $\rho_{\infty} = 1$.  Note that Eq. \eqref{eq:genalpha_flow} and Eq. \eqref{eq:genalpha_scalar} are completely decoupled.  We solve each of these using a two-stage predictor-multicorrector algorithm, as is common with Generalized-$\alpha$ methods.  For brevity we present below only the algorithm for Eq. \eqref{eq:genalpha_scalar}.\\

\noindent \textbf{Predictor Stage:} Set:
\\
\begin{equation}
{\dot{\bf{\Phi}}}_{n+1}^{(0)} = \frac{\gamma - 1}{\gamma} \dot{\bf{\Phi}}_n, \hspace{15pt} {\bf{{\Phi}}}_{n+1}^{(0)} = {\bf{{\Phi}}}_n, \hspace{15pt} \hat{{\bf{{\Phi}}}}_{n+1}^{(0)} = \hat{{\bf{{\Phi}}}}_n
\label{xpred}
\end{equation}
where:
\begin{equation}
\gamma = \frac{1}{2} +{\alpha_m} -{\alpha_f}
\end{equation}
\\

\noindent \textbf{Multicorrector Stage:} Repeat the following steps for $i = 1, 2, \ldots, i_{\text{max}}$:\\
\\
\noindent \textit{Step 1:} Evaluate iterates at the $\alpha$-levels:
\begin{align}
{\dot{\bf{\Phi}}}_{n+\alpha_m}^{(i)} &= {\dot{\bf{\Phi}}}_{n} + \alpha_m \left( \dot{\bf{\Phi}}_{n+1}^{(i-1)} - {\dot{\bf{\Phi}}}_{n} \right) \\
{\bf{\Phi}}_{n+\alpha_f}^{(i)} &= {\bf{\Phi}}_{n} + \alpha_f \left( {\bf{\Phi}}_{n+1}^{(i-1)} - {\bf{\Phi}}_{n} \right)
\end{align}

\noindent \textit{Step 2:} Use the solutions at the $\alpha$-levels to assemble the residual and the tangent matrix of the linear system:
{
\begin{displaymath}
\left[ \begin{array}{cc}
{\bf{A}}^{(i)} & {\bf{C}}^{(i)}  \\
{\bf{B}}^{(i)} & {\bf{D}}^{(i)}  \\
    \end{array} \right] \left[ \begin{array}{c}
\Delta \dot{\bf{\Phi}}_{n+1}^{(i)}\\
\Delta \hat{\bf{\Phi}}_{n+1}^{(i)}
\\ \end{array} \right] = -\left[ \begin{array}{c}
    {\bf{R}}^{(i)}_{\text{scalar,int}}\\
    {\bf{R}}^{(i)}_{\text{scalar,trace}} \\
 \end{array} \right]
\end{displaymath}
}
\\
where ${\bf{R}}^{(i)}_{\text{scalar,int}} = {\bf{R}}_{\text{scalar,int}}(\dot{{\bf{\Phi}}}^{(i)}_{n+\alpha_m},{\bf{\Phi}}^{(i)}_{n+\alpha_f},\hat{{\bf{\Phi}}}^{(i)}_{n+1},{\bf{U}}_n)$ and ${\bf{R}}^{(i)}_{\text{scalar,trace}} = {\bf{R}}_{\text{scalar,trace}}({\bf{\Phi}}^{(i)}_{n+\alpha_f},\hat{{\bf{\Phi}}}^{(i)}_{n+1},{\bf{U}}_n)$ are the residual equations associated with the turbulence model equation and the turbulence model conservativity condition, respectively, and:
\begin{align}
{\bf{A}}^{(i)} &= \alpha_m \frac{\partial {\bf{R}}_{\text{scalar,int}}(\dot{{\bf{\Phi}}}^{(i)}_{n+\alpha_m},{\bf{\Phi}}^{(i)}_{n+\alpha_f},\hat{{\bf{\Phi}}}^{(i)}_{n+1},{\bf{U}}_n)}{\partial \dot{\bf{\Phi}}_{n+\alpha_m}} + \alpha_f \gamma \Delta t_n \frac{\partial {\bf{R}}_{\text{scalar,int}}(\dot{{\bf{\Phi}}}^{(i)}_{n+\alpha_m},{\bf{\Phi}}^{(i)}_{n+\alpha_f},\hat{{\bf{\Phi}}}^{(i)}_{n+1},{\bf{U}}_n)}{\partial {\bf{\Phi}}_{n+\alpha_f}} \\
{\bf{B}}^{(i)} &= \alpha_f \gamma \Delta t_n \frac{\partial {\bf{R}}_{\text{scalar,trace}}({\bf{\Phi}}^{(i)}_{n+\alpha_f},\hat{{\bf{\Phi}}}^{(i)}_{n+1},{\bf{U}}_n)}{\partial {\bf{\Phi}}_{n+\alpha_f}} \\
{\bf{C}}^{(i)} &= \frac{\partial {\bf{R}}_{\text{scalar,int}}(\dot{{\bf{\Phi}}}^{(i)}_{n+\alpha_m},{\bf{\Phi}}^{(i)}_{n+\alpha_f},\hat{{\bf{\Phi}}}^{(i)}_{n+1},{\bf{U}}_n)}{\partial \hat{\bf{\Phi}}_{n+1}}\\
{\bf{D}}^{(i)} &= \frac{\partial {\bf{R}}_{\text{scalar,trace}}({\bf{\Phi}}^{(i)}_{n+\alpha_f},\hat{{\bf{\Phi}}}^{(i)}_{n+1},{\bf{U}}_n)}{\partial \hat{\bf{\Phi}}_{n+1}}
\end{align}
are the corresponding tangent matrices where $\Delta t_n = t_{n+1} - t_{n}$ is the time-step size.

\noindent \textit{Step 3:} Solve the linear system assembled in Step 2 for the update vectors $\Delta \dot{\bf{\Phi}}_{n+1}^{(i)}$ and $\Delta \hat{\bf{\Phi}}_{n+1}^{(i)}$, and update the iterates using the relations:
\begin{align}
{\dot{\bf{\Phi}}}_{n+1}^{(i)} &= {\dot{\bf{\Phi}}}_{n+1}^{(i-1)} +  \Delta \dot{\bf{\Phi}}_{n+1}^{(i)}\\
{\bf{\Phi}}_{n+1}^{(i)} &= {\bf{\Phi}}_{n+1}^{(i-1)} +  \gamma \Delta t_n \Delta \dot{\bf{\Phi}}_{n+1}^{(i)}\\
{\bf{\hat{\Phi}}}_{n+1}^{(i)} &= {\bf{\hat{\Phi}}}_{n+1}^{(i-1)} +  \Delta \hat{\bf{\Phi}}_{n+1}^{(i)}
\end{align}

\noindent \textit{Step 4:} If the norms of both ${\bf{R}}^{(i)}_{\text{scalar,int}}$ and ${\bf{R}}^{(i)}_{\text{scalar,trace}}$ are less than some prescribed tolerance (e.g., $10^{-6}$ of their initial values), terminate the iterative loop and set ${\dot{\bf{\Phi}}}_{n+1} = {\dot{\bf{\Phi}}}_{n+1}^{(i)}$, ${{\bf{\Phi}}}_{n+1} = {{\bf{\Phi}}}_{n+1}^{(i)}$, and ${\hat{\bf{\Phi}}}_{n+1} = {\hat{\bf{\Phi}}}_{n+1}^{(i)}$.

\section{Static Condensation}

With each iterate of the predictor-multicorrector algorithm presented in the previous section, we must solve a linear system of the form:
{
\begin{displaymath}
\left[ \begin{array}{cc}
{\bf{A}}^{(i)} & {\bf{C}}^{(i)}  \\
{\bf{B}}^{(i)} & {\bf{D}}^{(i)}  \\
    \end{array} \right] \left[ \begin{array}{c}
\Delta \dot{\bf{\Phi}}_{n+1}^{(i)}\\
\Delta \hat{\bf{\Phi}}_{n+1}^{(i)}
\\ \end{array} \right] = -\left[ \begin{array}{c}
    {\bf{R}}^{(i)}_{\text{scalar,int}}\\
    {\bf{R}}^{(i)}_{\text{scalar,trace}} \\
 \end{array} \right]
\end{displaymath}
}
\\
It should be noted that the matrix ${\bf{A}}^{(i)}$ appearing in the above system is block diagonal.  Consequently, we can use static condensation to significantly reduce the computational cost associated with solving the system, a particular advantage of the HDG approach over a classical DG method.  Namely, we can first solve the matrix system:
\begin{align*}
&{\bf{S}}^{(i)} \Delta {\hat{\bf{\Phi}}}^{(i)}_{n+1} = - {\bf{R}}^{(i)}_{\text{scalar,trace}} + {\bf{B}}^{(i)} \left( {\bf{A}}^{(i)} \right)^{-1} {\bf{R}}^{(i)}_{\text{scalar,int}}
\end{align*}
for the trace degrees-of-freedom $\Delta {\hat{\bf{\Phi}}}^{(i)}_{n+1}$ where:
\begin{align*}
&{\bf{S}}^{(i)} = {\bf{D}}^{(i)}-{\bf{B}}^{(i)} \left( {\bf{A}}^{(i)} \right)^{-1} {\bf{C}}^{(i)}
\end{align*}
is the Schur complement, and then we can solve the matrix system:
\begin{align*}
{\bf{A}}^{(i)} \Delta {\dot{\bf{\Phi}}}^{(i)}_{n+1} = - {\bf{R}}^{(i)}_{\text{scalar,int}} - {\bf{C}}^{(i)} \Delta {\hat{\bf{\Phi}}}^{(i)}_{n+1}
\end{align*}
for the interior degrees-of-freedom $\Delta {\dot{\bf{\Phi}}}^{(i)}_{n+1}$.  Since ${\bf{A}}^{(i)}$ is block diagonal, we can construct the Schur complement ${\bf{S}}^{(i)}$ quite efficiently.  In fact, the Schur complement can be formed and assembled element-wise in a standard element assembly routine \cite{cgorhdg}.  It should also be noted that the interior degrees-of-freedom $\Delta {\dot{\bf{\Phi}}}^{(i)}_{n+1}$ may be attained via local element-wise solves after the trace degrees-of-freedom $\Delta {\hat{\bf{\Phi}}}^{(i)}_{n+1}$ have been solved for as ${\bf{A}}^{(i)}$ is block diagonal.

\section{Numerical Results}

Now that we have presented our fully-discrete HDG formulation for the incompressible RANS equations coupled with the Spalart-Allmaras one equation turbulence model, we now perform verification of our method using a selection of two-dimensional example problems.  We first analyze the accuracy of our formulation using the method of manufactured solutions.  We then analyze the effectiveness of our formulation using three common benchmark problems: flow in a turbulent channel, flow over a backward facing step, and flow over a NACA 0012 airfoil at moderate Reynolds number and angle of attack.  Each of the example problems we consider are steady in the sense that the mean velocity and pressure fields are independent of time.  Nonetheless, to solve the example problems, we employed a transient analysis until a steady state flow solution was attained.  All simulations were initialized using a steady state Stokes flow solution for the mean flow field and a steady state diffusion solution for the working viscosity.  We also considered alternative initial conditions and found that the final steady state flow solution was independent of the initial condition in all cases.\\


\subsection{Manufactured Solution}

As a first numerical experiment, we consider a two-dimensional manufactured vortex solution for the mean flow field and a manufactured solution for the working viscosity.  Namely, the flow domain for this experiment is:
\[
\Omega = [0,1]^2
\]
the mean velocity field is:
\begin{equation}
\overline{{\bf{u}}} = \left[
\begin{aligned}
&-2 x^2 e^x (-y^2+y) (2y-1) (x-1)^2\\
&-x y^2 e^x (x (x+3)-2) (x-1) (y-1)^2 
\end{aligned}
\right]
\end{equation} 
the mean pressure field is:
\begin{equation}
\begin{aligned}
\overline{p} = \text{sin}(\pi x) \text{sin}(\pi y)
\end{aligned}
\end{equation}
and the working viscosity is:
\begin{equation}
\begin{aligned}
\tilde{\nu} = \theta\text{sin}(\pi x) \text{sin}(\pi y)
\end{aligned}
\end{equation}
where $\theta$ is a constant to be specified.  The density $\rho$ is set to 1, and the viscosity $\nu$ is left unspecified.
\begin{figure}[t!]
	\begin{subfigure}{0.33\linewidth}
  		\centering \includegraphics[scale=0.81]{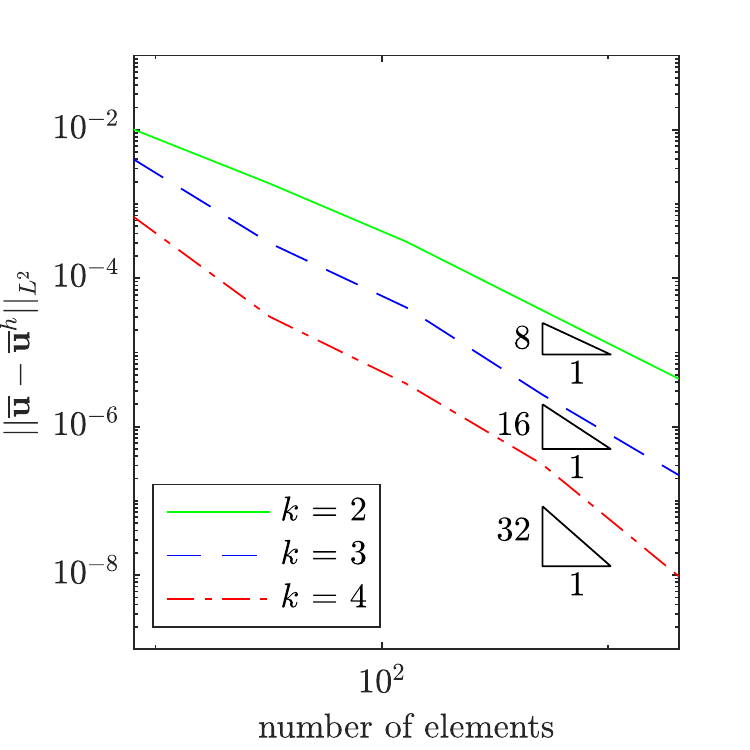}
  		\caption{Mean Velocity}
	\end{subfigure}
		\begin{subfigure}{0.33\linewidth}
  		\centering \includegraphics[scale=0.81]{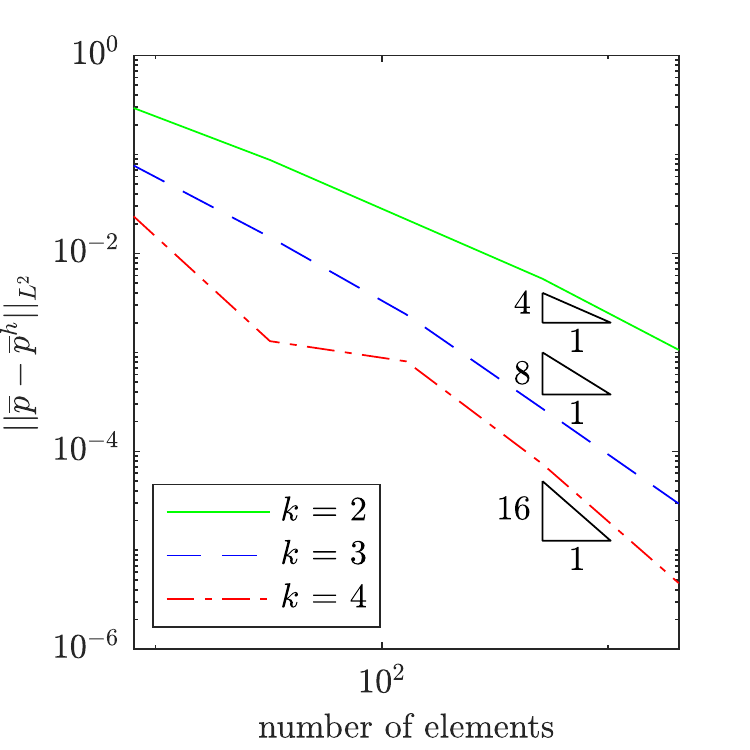}
  		\caption{Mean Pressure}
	\end{subfigure}
	\begin{subfigure}{0.33\linewidth}
		\centering
  		\includegraphics[scale=0.81]{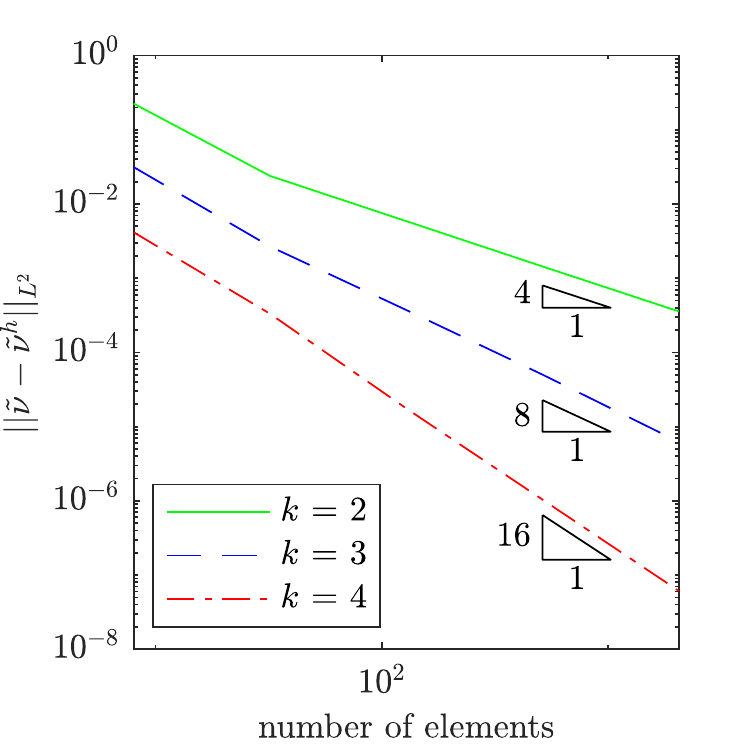}
  		\caption{Working Viscosity}
	\end{subfigure}
	\caption{Convergence for the method of manufactured solutions problem with $\nu =$ 1e-2, $\theta =$ 1e0, a mean velocity polynomial degree of $k$, a mean pressure polynomial degree of $k -1$, and a working viscosity polynomial degree of $k - 1$ for $k = 2, 3, 4$.}
	\label{convergence1}
\end{figure}
The corresponding forcing in the momentum equation is:
\begin{equation}
{\bf{f}} = \nabla \cdot (\overline{{\bf{u}}} \otimes \overline{{\bf{u}}}) + \nabla \overline{p} - \nabla \cdot (2(\nu+\nu_T(\tilde{\nu}))\nabla^S \overline{{\bf{u}}}) 
\end{equation}
and the forcing in the Spalart-Allmaras turbulence model equation is:
\begin{equation}
\begin{aligned}
f = \overline{{\bf{u}}} \cdot \nabla \tilde{\nu} - \nabla \cdot \frac{1}{\sigma} \left[ (\nu + \tilde{\nu} ) \nabla \tilde{\nu} \right] - C_{b1} \tilde{S} \tilde{\nu} + C_{w1} f_{w} \frac{\tilde{\nu}^2}{d^2} - \frac{C_{b2}}{\sigma} \nabla \tilde{\nu} \cdot \nabla \tilde{\nu}
\end{aligned}
\end{equation}
Homogeneous Dirichlet boundary conditions are applied along the boundary $\partial \Omega$ for both the mean velocity field and the working viscosity, and the condition $\int_{\Omega} \overline{p} d\Omega = 0$ is enforced.\\
\\
\noindent To assess the accuracy of our method, we have solved the manufactured solution problem for a selection of values for $\theta$ and $\nu$ as well as a series of meshes and several different polynomial degrees for the mean flow field and the working viscosity.  In Fig. \ref{convergence1}, results are displayed for $\nu =$ 1e-2, $\theta =$ 1e0, a mean velocity polynomial degree of $k$, a mean pressure polynomial degree of $k - 1$, and a working viscosity polynomial degree of $k - 1$ for $k = 2, 3, 4$.  From the plots, it is apparent that the mean velocity field, the mean pressure field, and the working viscosity converge at the expected rates of $k+1$, $k$, and $k$ in the $L^2$-norm.  We observed this same behavior for other values of $\nu$ and $\theta$ as well.

\begin{figure}[t!]
	\begin{subfigure}{0.33\linewidth}
  		\centering \includegraphics[scale=0.81]{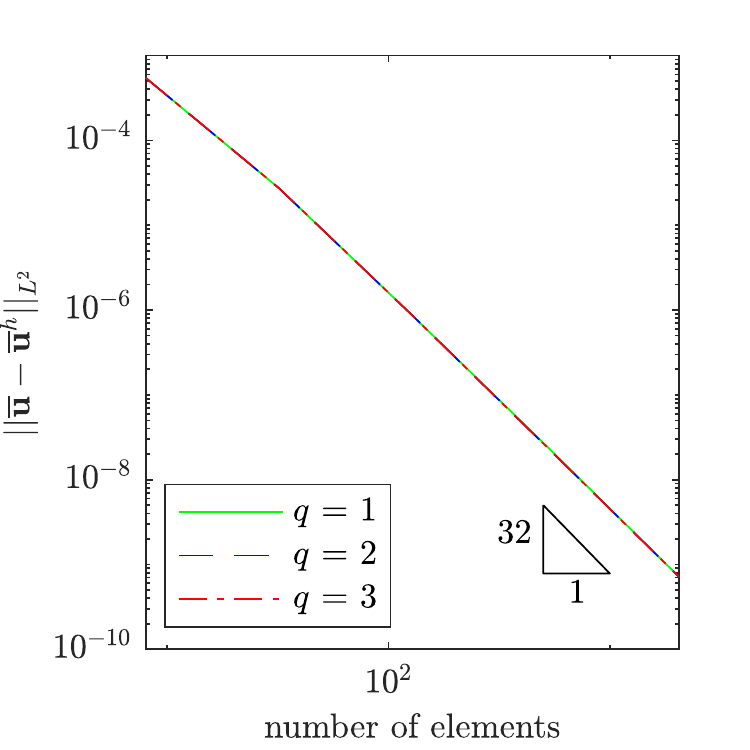}
  		\caption{Mean Velocity}
	\end{subfigure}
		\begin{subfigure}{0.33\linewidth}
  		\centering \includegraphics[scale=0.81]{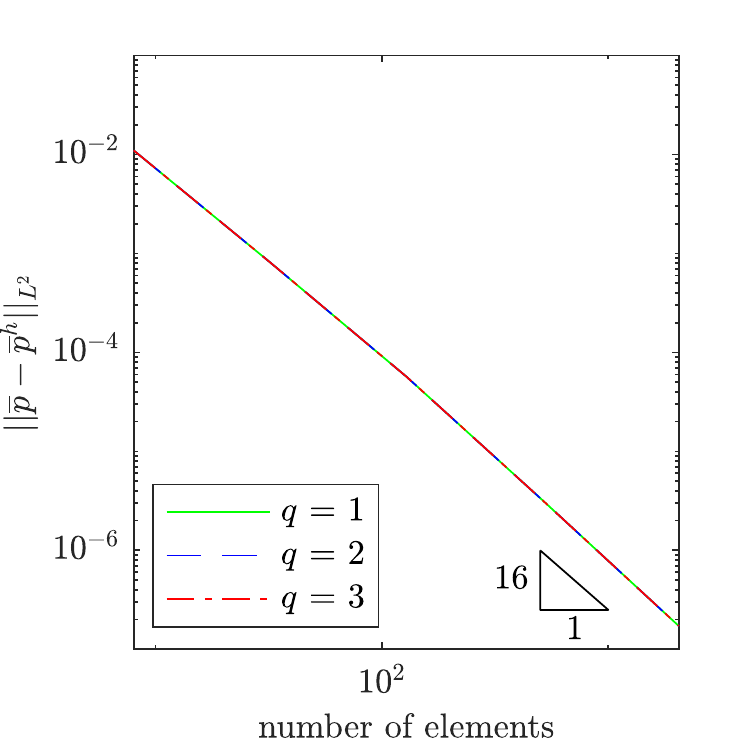}
  		\caption{Mean Pressure}
	\end{subfigure}
	\begin{subfigure}{0.33\linewidth}
		\centering
  		\includegraphics[scale=0.81]{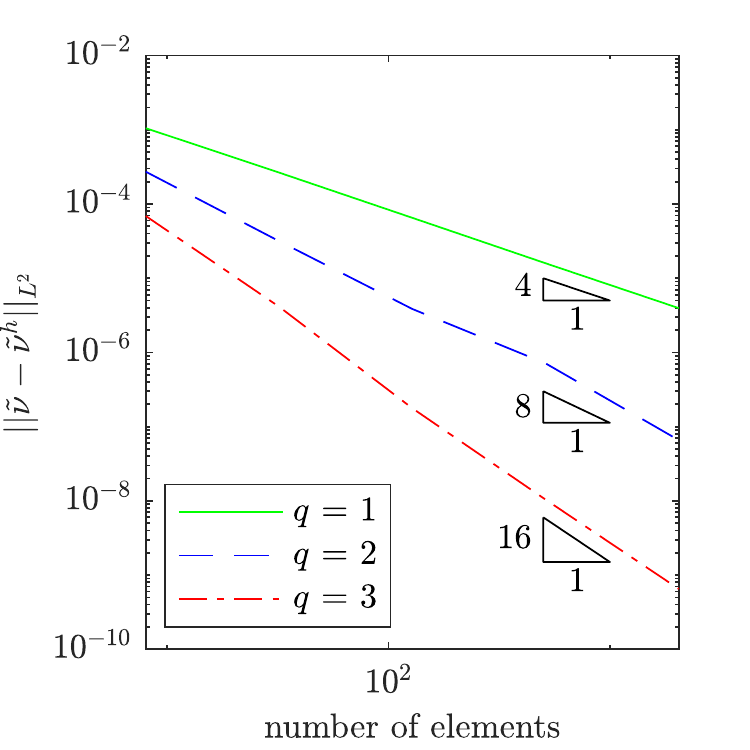}
  		\caption{Working Viscosity}
	\end{subfigure}
	\caption{Convergence for the method of manufactured solutions problem with $\nu =$ 1e-2, $\theta =$ 1e-2, a mean velocity polynomial degree of 4, a mean pressure polynomial degree of 3, and a working viscosity polynomial degree of $q$ for $q = 1, 2, 3$.}
	\label{convergence2}
\end{figure}

\begin{figure}[t!]
	\begin{subfigure}{0.33\linewidth}
  		\centering \includegraphics[scale=0.81]{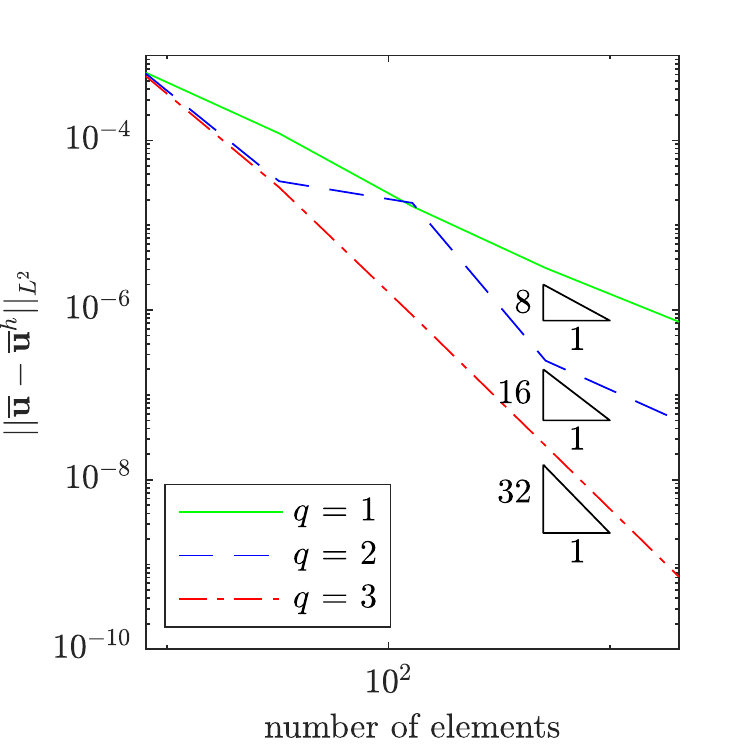}
  		\caption{Mean Velocity}
	\end{subfigure}
		\begin{subfigure}{0.33\linewidth}
  		\centering \includegraphics[scale=0.81]{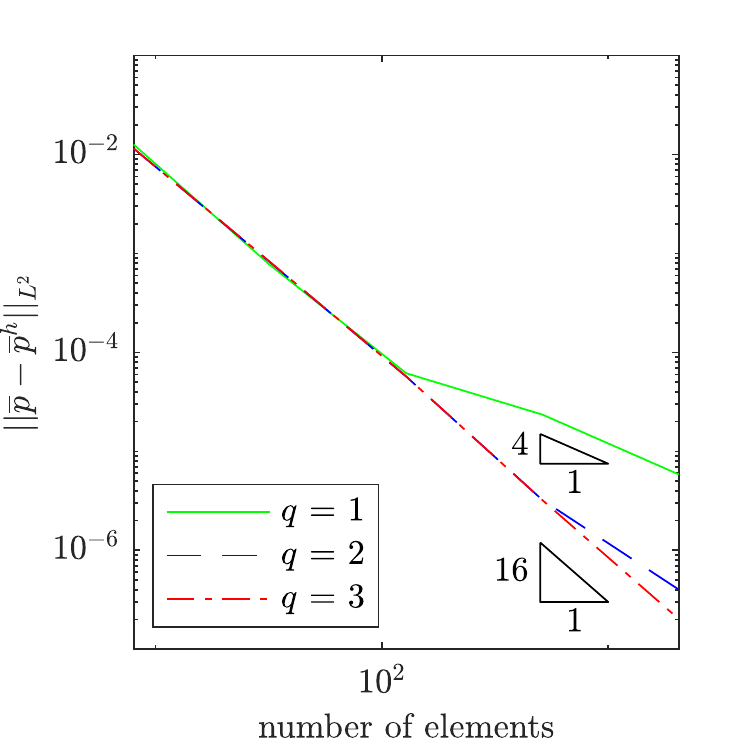}
  		\caption{Mean Pressure}
	\end{subfigure}
	\begin{subfigure}{0.33\linewidth}
		\centering
  		\includegraphics[scale=0.81]{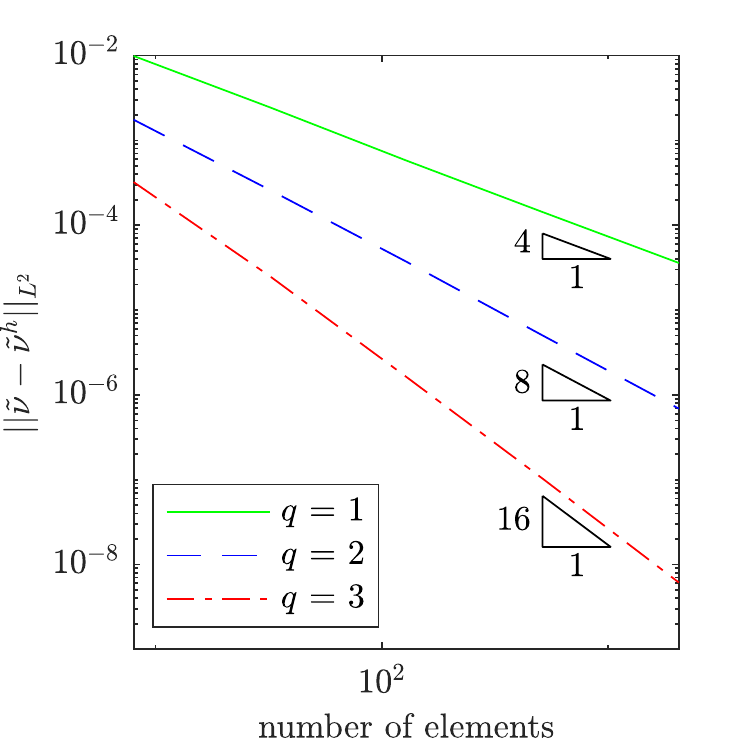}
  		\caption{Working Variable}
	\end{subfigure}
	\caption{Convergence for the method of manufactured solutions problem with $\nu =$ 1e-2, $\theta =$ 1e-1, a mean velocity polynomial degree of 4, a mean pressure polynomial degree of 3, and a working variable polynomial degree of $q$ for $q = 1, 2, 3$.}
	\label{convergence3}
\end{figure}

\noindent To assess the impact of the polynomial degree of the working viscosity on the accuracy of the mean flow field, we also examined the impact of using a lower polynomial degree for the working viscosity than for the mean flow field.  While theoretically, we expect the rate of convergence of the mean flow field to depend on the polynomial degree of both the mean flow field and the working viscosity, for certain configurations, we hypothesize that a lower polynomial degree can be employed for the working viscosity than for the mean flow field.  In Figs. \ref{convergence2} and \ref{convergence3}, results are displayed for $\nu =$ 1e-2, $\theta =$ 1e-2 and 1e-1 respectively, a mean velocity polynomial degree of 4, a mean pressure polynomial degree of 3, and a working viscosity polynomial degree of $q$ for $q = 1, 2, 3$.  From the plots, we observe that optimal convergence rates are attained for the mean velocity and pressure fields for $\theta =$ 1e-2, but reduced convergence rates are attained for $\theta =$ 1e-1 for higher levels of mesh refinement.  These plots and further studies suggest that optimal convergence rates for the mean velocity and pressure fields are attained only pre-asymptotically, with a pre-asymptotic range whose size decreases with increasing $\theta$, unless the polynomial degree for the working viscosity is chosen to be at most one order less than the polynomial degree associated with the mean velocity.  Nonetheless, for many problems of interest, including the common benchmark problems considered here, we have found that accurate mean velocity and pressure results are attainable using a low polynomial degree for the working viscosity.\\


\newpage
\subsection{Turbulent Channel Flow at $Re_\tau = $ 550}
The next example we consider is turbulent flow in a channel at $\text{Re}_{{\tau}}$ = 550, where $\text{Re}_{{\tau}}$ is the Reynolds number based on the friction velocity  and the channel half width.  The domain is rectangular with dimensions 0.5 x 1.0 in the streamwise and wall normal directions, respectively.  Only half the channel is modeled.  For the mean flow field, periodic boundary conditions are applied in the streamwise direction, no-slip Dirichlet boundary conditions are imposed at the wall, and a symmetry condition is applied at the middle of the channel.  For the working viscosity, periodic boundary conditions are applied in the streamwise direction, a homogeneous Dirichlet boundary condition is applied at the wall, and a homogeneous Neumann boundary condition is applied at the middle of the channel.    32 elements are employed in the wall-normal direction, and non-uniform grid spacing is utilized so that the boundary layer can be resolved.  The polynomial degrees of the mean velocity, mean pressure, and working viscosity are taken to be 3, 2, and 1, respectively.  An externally imposed pressure gradient in the streamwise direction is imposed, and the value of this gradient is selected to ensure the correct $\text{Re}_{{\tau}}$ is achieved.  The kinematic viscosity is selected to be $\nu$ = 1e-4.\\
\begin{figure}[t!]
	\centering
	\includegraphics[scale=0.8]{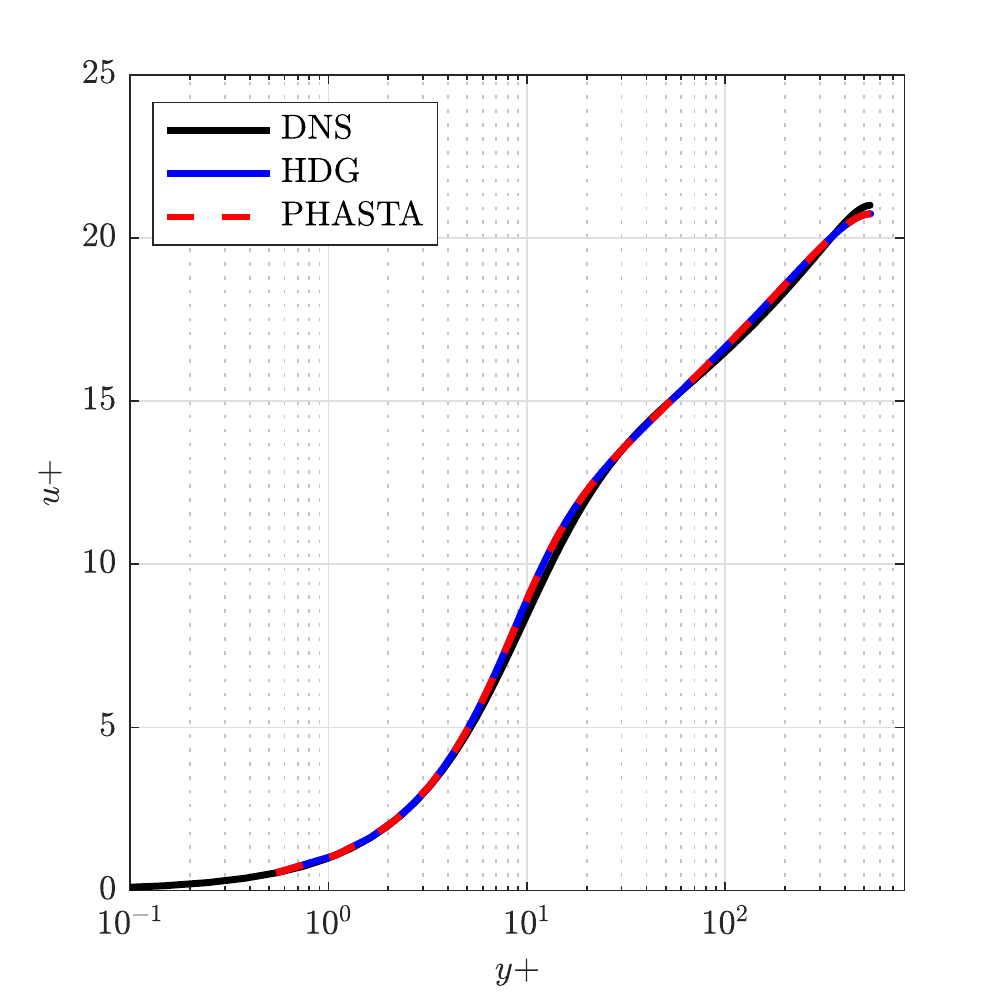}
	\caption{Comparison of the velocity profile for DNS and RANS with the Spalart-Allmaras model for a turbulent channel flow at $\text{Re}_{\tau}$ = 550.}
	\label{wall}
\end{figure}
\\
\noindent In Fig. \ref{wall}, the mean velocity in the streamwise direction normalized by the friction velocity $u^+$ versus the non-dimensional distance to the wall in wall-units $y^+$ is displayed.  Results attained using the HDG method presented in this paper are compared with those attained using a stabilized finite element implementation of the Spalart-Allmaras model within the open source CFD library PHASTA \cite{phasta1,phasta2}, as well as the DNS data of Lee and Moser \cite{moser_direct_2015}.  A grid refinement study was conducted to ensure the stabilized finite element results were fully converged.  The HDG results and the stabilized finite element results are nearly indistinguishable, which is expected as they are solving the same set of partial differential equations.  The HDG results also closely approximate the DNS results, which is a testament of the power of the Spalart-Allmaras model more-so than the discretization since the DNS solves a different set of partial differential equations.

\subsection{Turbulent Backward Facing Step at $Re =$ 36,000}
\noindent We next consider turbulent flow over a backward facing step.  This problem is challenging because of the separation region that occurs downstream of the step at sufficiently large enough Reynolds numbers.  The reattachment length in particular is notoriously difficult to predict.  It is well-known that the Spalart-Allmaras model is not ideal for this problem, as it was created for attached wall bounded flows.  However, our goal in exploring the turbulent backward facing step is not to see how well our method is able to match DNS data but rather how well it is able to match Spalart-Allmaras data in the literature and how quickly it converges with mesh refinement and polynomial degree elevation.\\
\\
\noindent For all of the results reported here, a Reynolds number of $Re =$ 36,000 based on the inflow bulk velocity and step height is employed.  Moreover, the kinematic viscosity is selected to be $\nu =$ 2.77777e-5, and the incoming bulk velocity is selected to be $1$.  The remainder of the problem setup is as described on the NASA turbulence modeling website \cite{nasa}.  Dirichlet boundary conditions for both the mean velocity and working viscosity are applied at the inlet, with the turbulent viscosity set to be four times the kinematic viscosity.  Homogeneous no-slip boundary conditions for the mean velocity and a homogeneous Dirichlet boundary condition for the working viscosity are set at the walls.  Zero-traction boundary conditions for the mean velocity and a homogeneous Neumann boundary condition for the working viscosity are prescribed at the outlet.\\
\begin{figure}[t!]
	\centering
	\includegraphics[scale=.3]{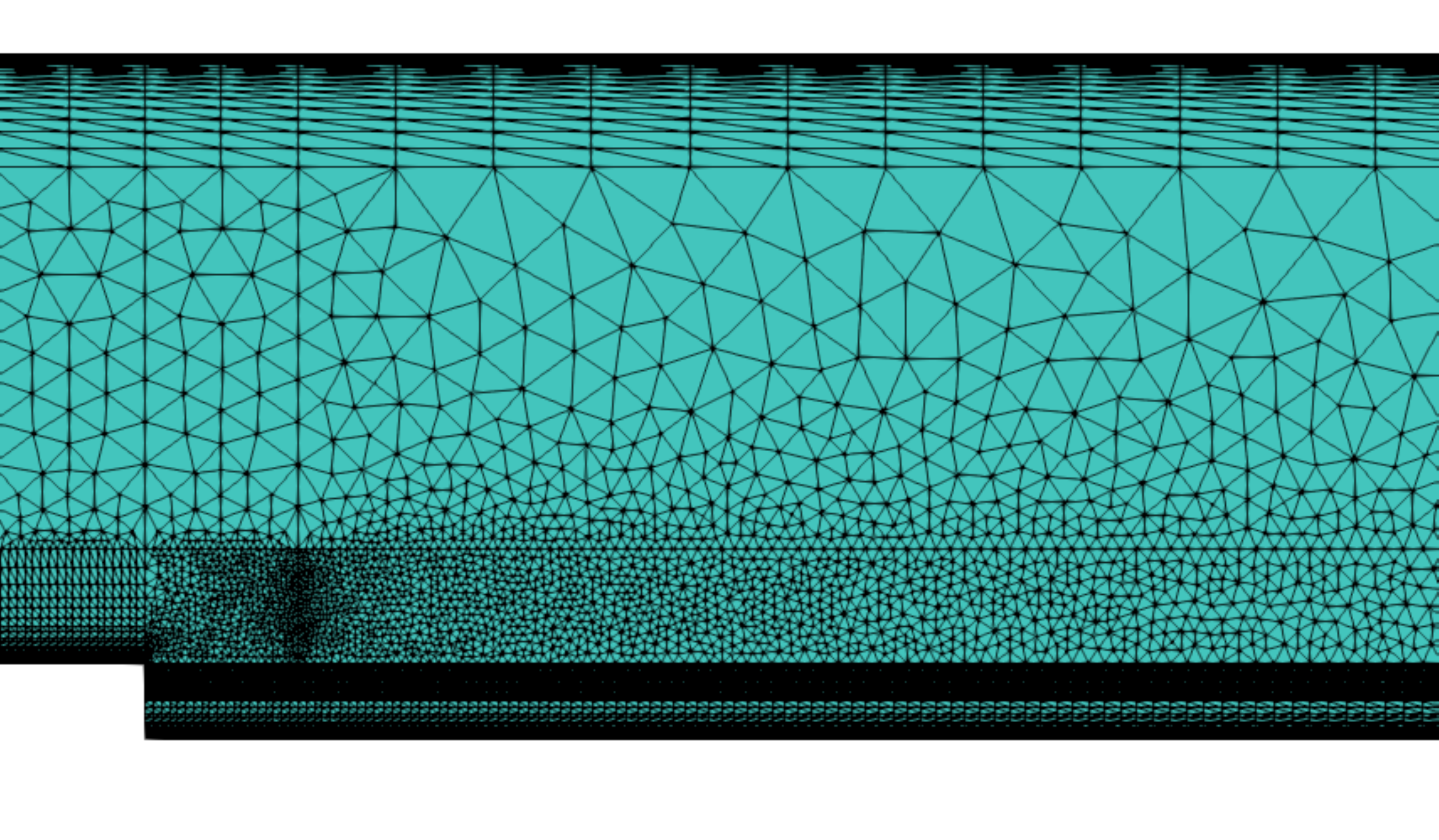}
	\caption{Computational mesh associated with first results for simulating turbulent flow over a backward facing step.}
	\label{backmesh}
\end{figure}
\begin{figure}[t!]
 	\centering
	\begin{subfigure}{0.49\linewidth}
		\centering
  		\includegraphics[scale=.43]{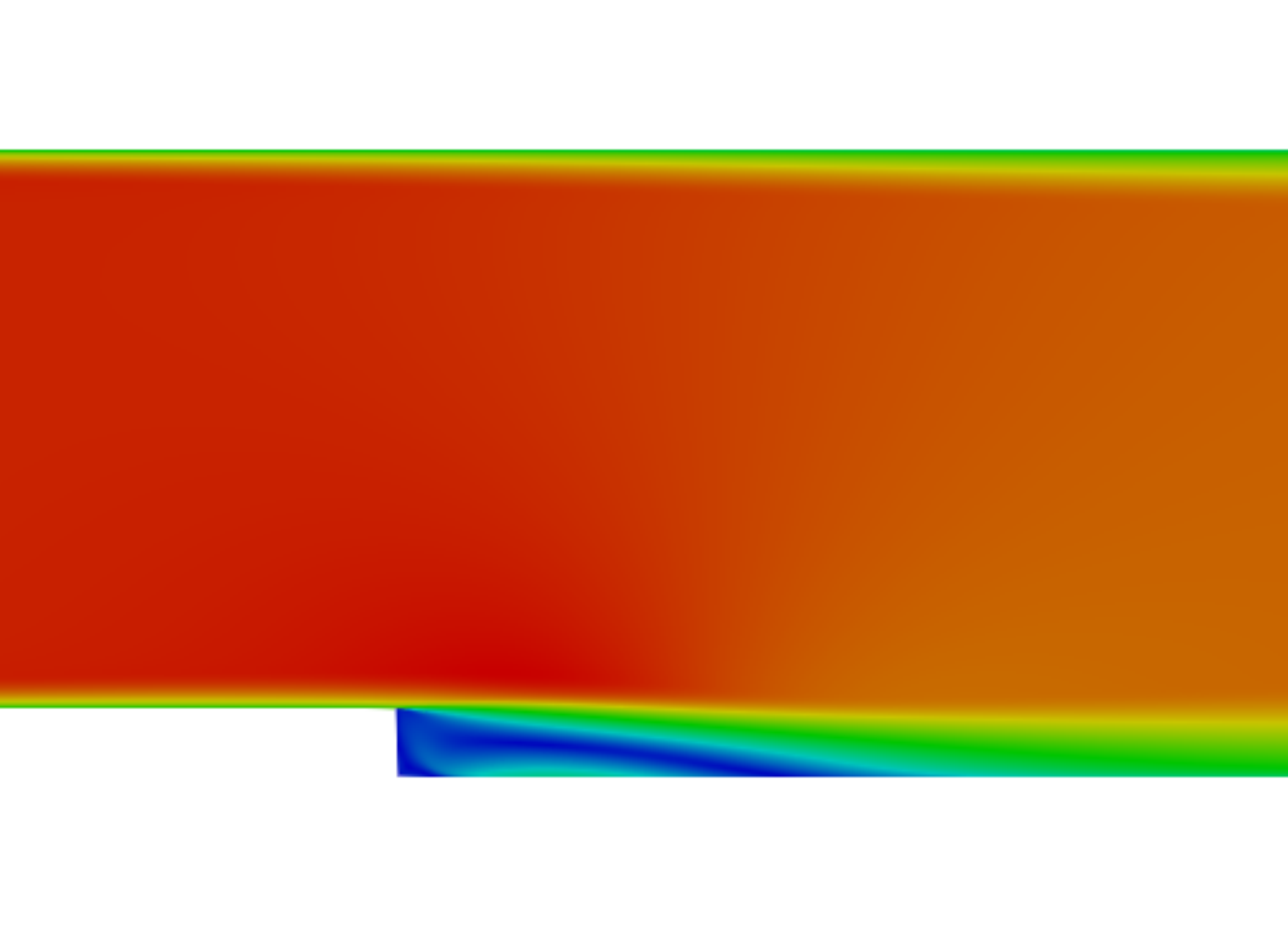}
  		\caption{Magnitude of Mean Velocity $\overline{{\bf{u}}}$}
	\end{subfigure}
		\begin{subfigure}{0.49\linewidth}
		\centering
  		\includegraphics[scale=.43]{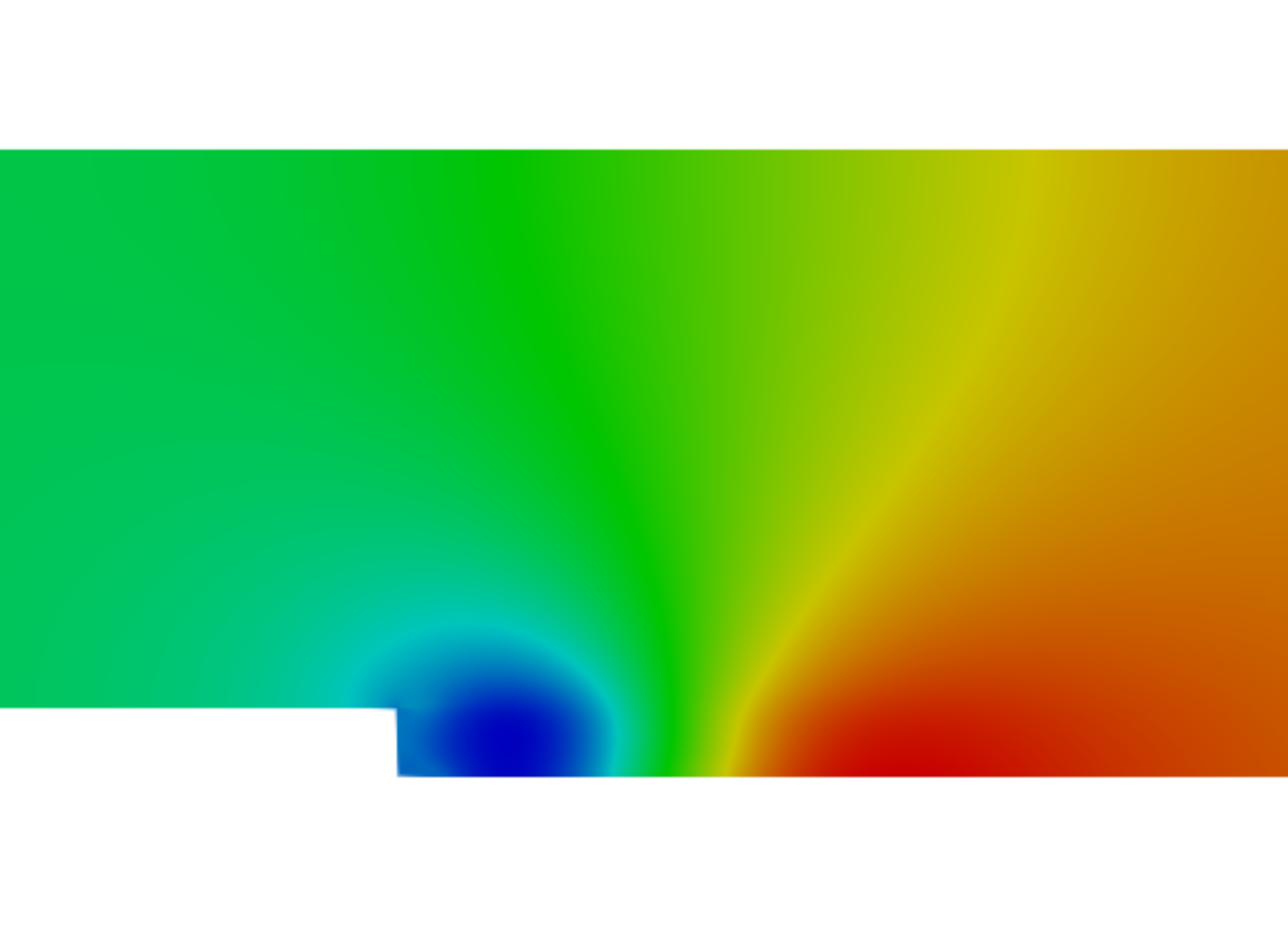}
  		\caption{Mean Pressure $\overline{p}$}
	\end{subfigure}
	\begin{subfigure}{1.0\linewidth}
		\centering
  		\includegraphics[scale=.43]{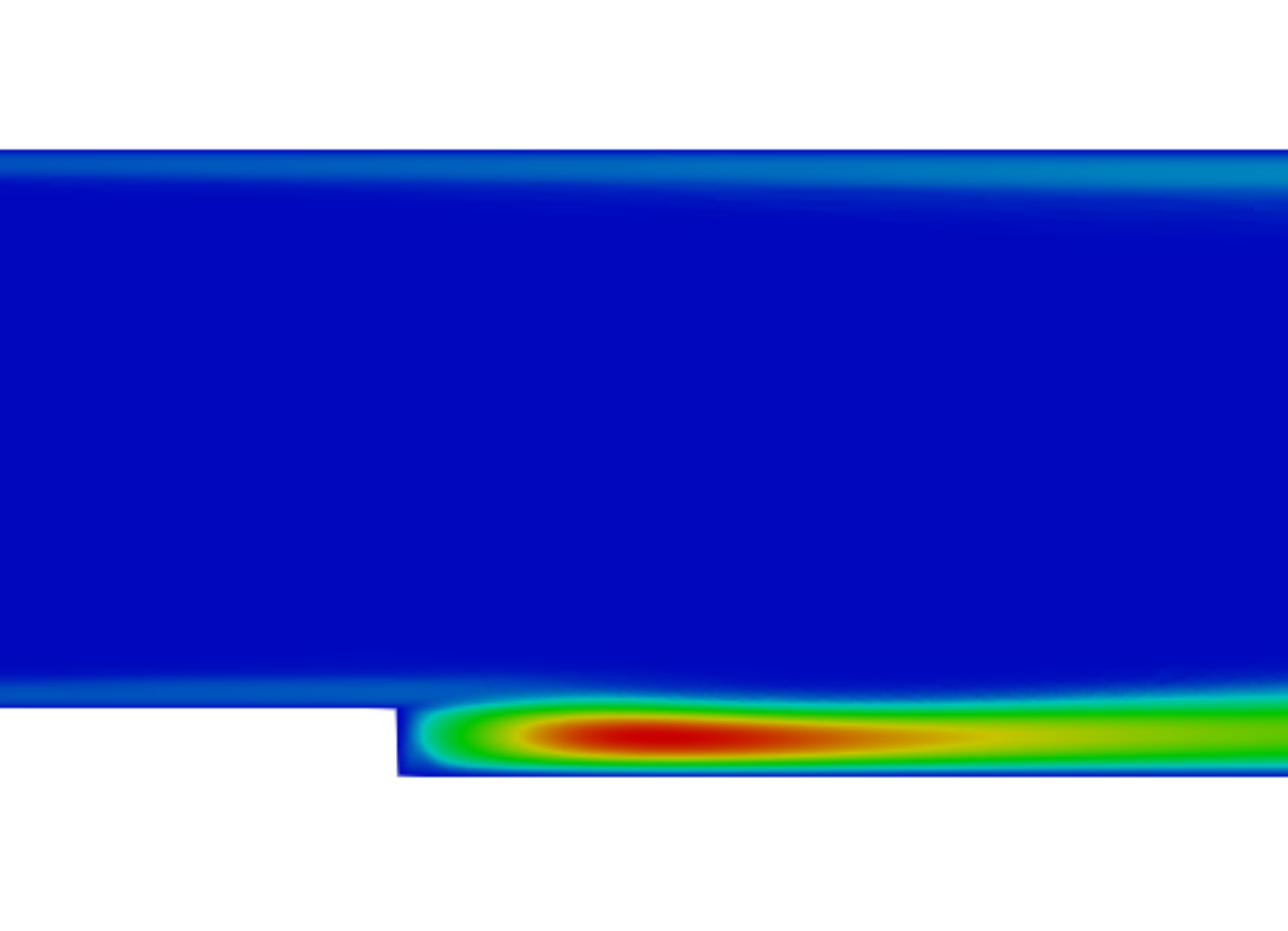}
  		\caption{Normalized Eddy Viscosity $\frac{\nu_T}{\nu}$}
	\end{subfigure}
	\caption{Computed velocity, pressure, and normalized eddy viscosity contours for flow over the backwards facing step at $Re = 36,000$ using mean velocity, mean pressure, and working viscosity polynomial degrees of 3, 2, and 1 respectively and the mesh displayed in Fig. \ref{backmesh}.}
	\label{backstep36}
\end{figure}
\\
\begin{figure}[t!]
    \centering
    \begin{subfigure}[b]{0.49\textwidth}
        \includegraphics[scale  = 1]{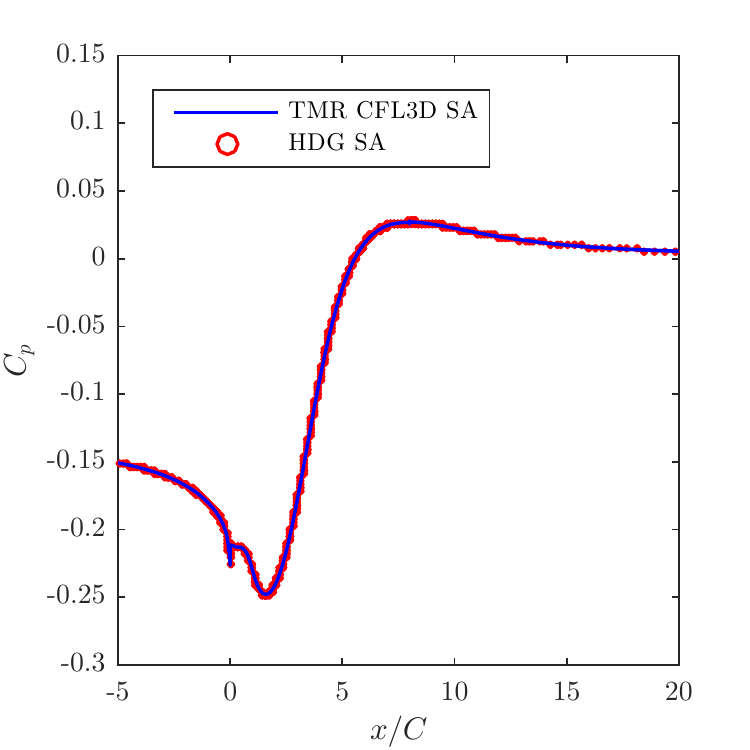}
        \caption{Surface Pressure Coefficient}
    \end{subfigure}
    \begin{subfigure}[b]{.49\textwidth}
        \includegraphics[scale = 1]{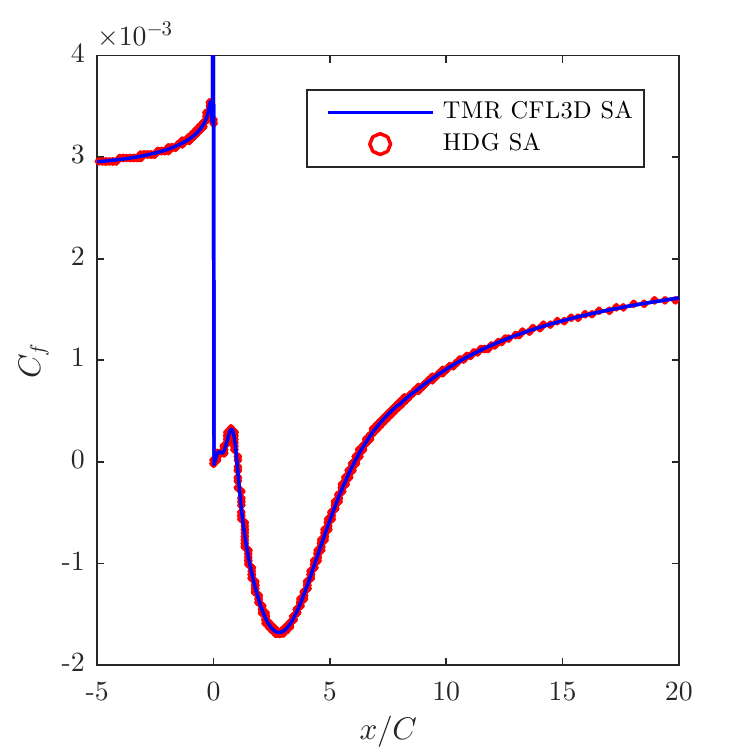}
        \caption{Surface Skin Friction Coefficient}
    \end{subfigure}
    \caption{Computed surface pressure and skin friction coefficients for flow over the backwards facing step at $Re = 36,000$ using mean velocity, mean pressure, and working viscosity polynomial degrees of 3, 2, and 1 respectively and the mesh displayed in Fig. \ref{backmesh}.  HDG results are plotted with finite volume results obtained using CFL3D and a mesh of 66,049 elements.}
        \label{profiles}
\end{figure}
\\
\begin{figure}[t!]
	\centering
	\includegraphics[scale=1.0]{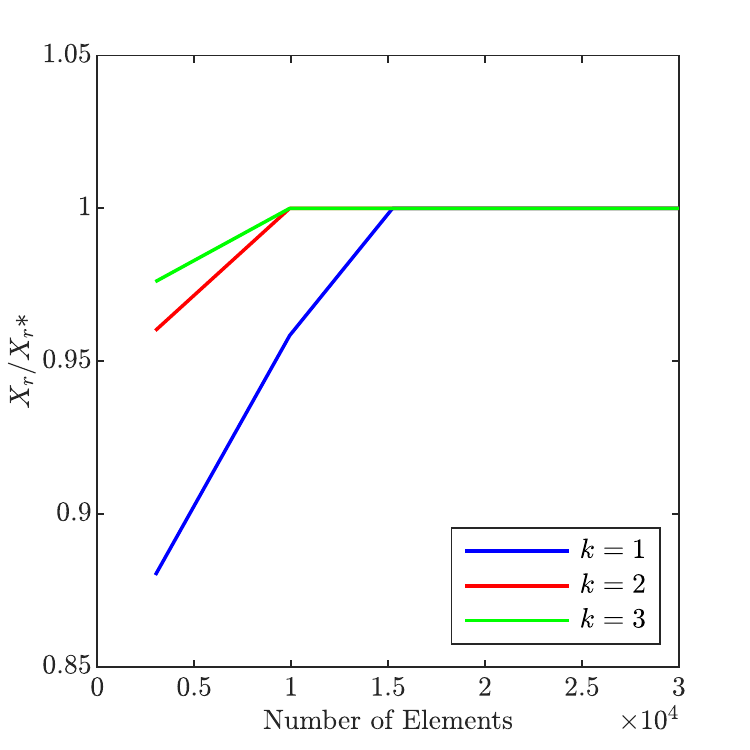}
	\caption{Refinement study for normalized reattachment length for flow over the backwards facing step at $Re = 36,000$.  In the figure, $X_r$ denotes the reattachment length for a given mesh size and polynomial degree while $X_{r^*}$ denotes the converged reattachment length of approximately 6.1.}
	\label{backmesh2c}
\end{figure}
\\
\noindent Our first results correspond to a computational mesh consisting of 15,325 elements and mean velocity, mean pressure, and working viscosity polynomial degrees of 3, 2, and 1, respectively.  The computational mesh is displayed in Fig. \ref{backmesh}, and attained steady state solutions for the mean velocity, mean pressure, and normalized eddy viscosity are displayed in Fig. \ref{backstep36}.  In order to assess the accuracy of our results, we compare in Fig. \ref{profiles} the obtained surface pressure and skin friction coefficients on the bottom surface of the domain with results publicly available on the NASA turbulence modeling website which were obtained using the NASA code CFL3D.  The CFL3D results were obtained using a second-order finite volume method and a mesh of 66,049 elements, while only 15,325 elements were employed for the HDG results.  Despite this discrepancy in resolution, the HDG results are visually indistinguishable from the CFL3D results.  Most notably, both the HDG results and the CFL3D results predict that reattachment occurs at $x \approx 6.1$.\\
\\
\noindent We next examine the impact of mesh refinement and polynomial degree elevation on the accuracy of our HDG method with a particular focus on reattachment length.  In Fig. \ref{backmesh2c}, the predicted reattachment length is reported for a series of refined meshes and for mean velocity, mean pressure, and working viscosity polynomial degrees of $k$, $k -1$, and 1 for $k = 1, 2, 3$.  Note that the predicted reattachment length quickly converges with mesh refinement for each of the considered polynomial degrees.  Further note that the accuracy is significantly improved with polynomial degree elevation.

\subsection{Turbulent Flow Over a NACA 0012 Airfoil at $\alpha$ = $10^{\circ}$}

\begin{figure}[t!]
	\centering
	    \begin{subfigure}[b]{1\textwidth}
	    \centering
        			\includegraphics[scale  = .35]{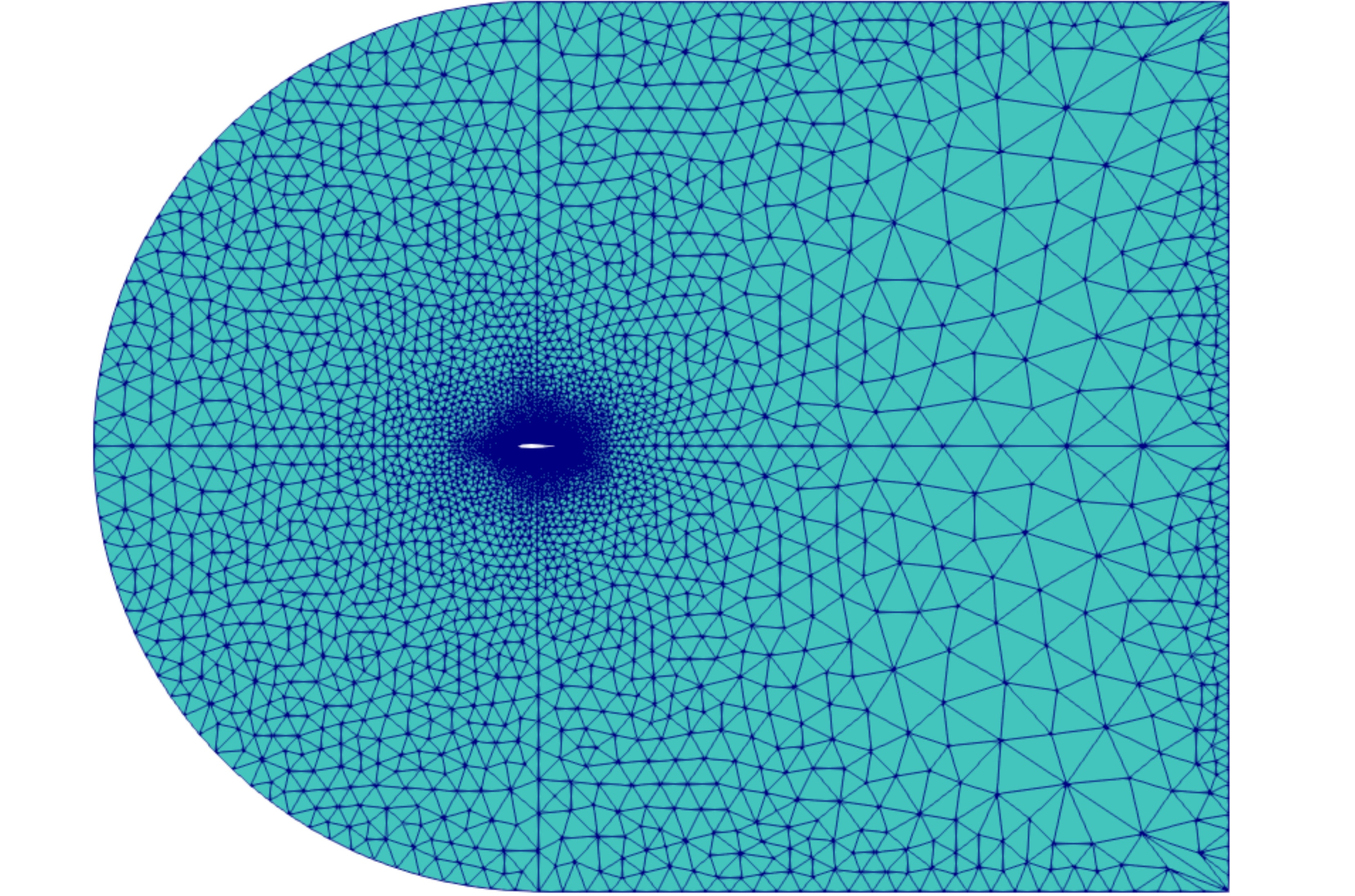}
       			 \caption{Far-Field View of Mesh}
   	    \end{subfigure}
	    	    \begin{subfigure}[b]{1\textwidth}
		    \centering
        			\includegraphics[scale  = .35]{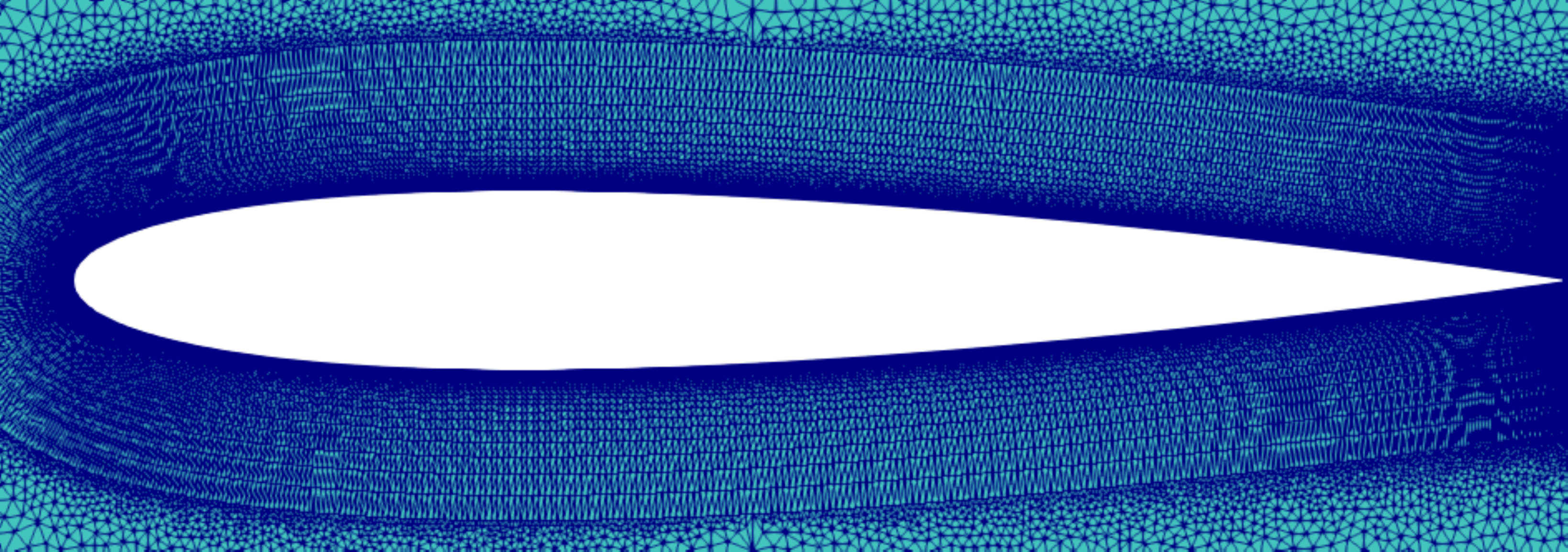}
       			 \caption{Boundary Layer View of Mesh}
   	    \end{subfigure}
	    	    	    \begin{subfigure}[b]{1\textwidth}
		    \centering
        			\includegraphics[scale  = .35]{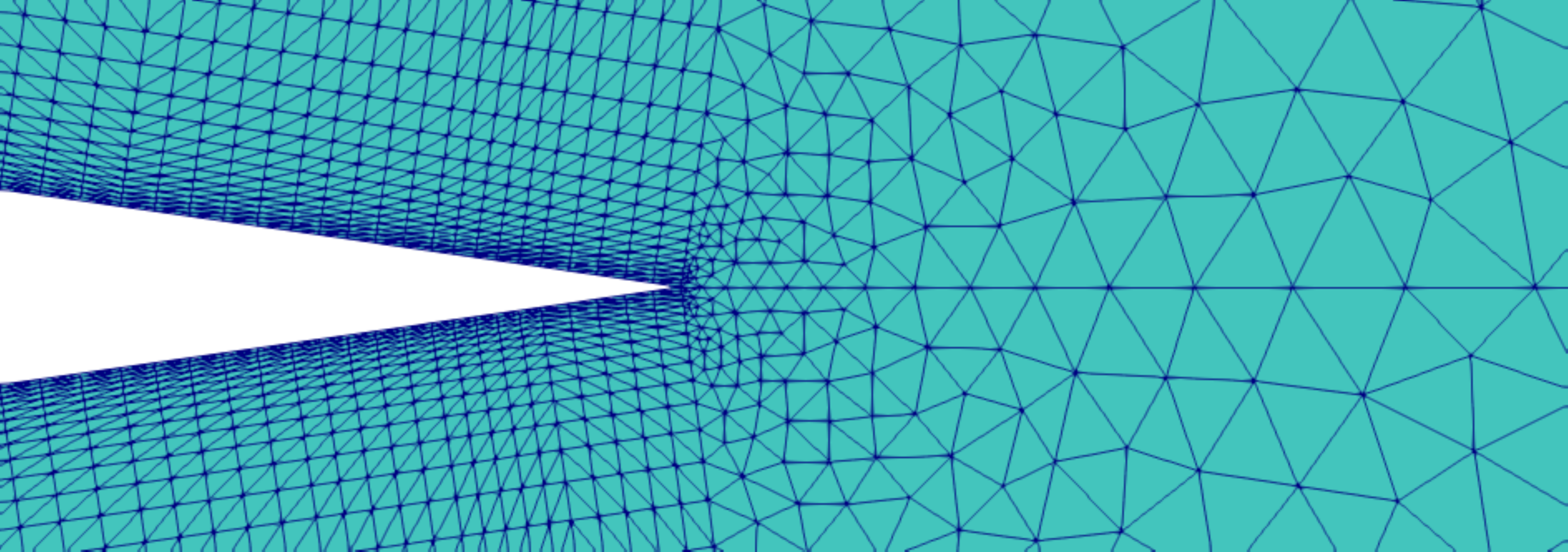}
       			 \caption{Zoomed In Trailing Edge View of Mesh Highlighting Boundary Layer Stretching}
   	    \end{subfigure}
	    \caption{Computational mesh associated with results for simulating turbulent flow over a NACA 0012 airfoil at $Re = 3,000,000$ and $\alpha=10^{\circ}$.}
	\label{nacamesh}
\end{figure}

\begin{figure}[t!]
 	\centering
	\begin{subfigure}{0.49\linewidth}
		\centering
  		\includegraphics[scale=0.85]{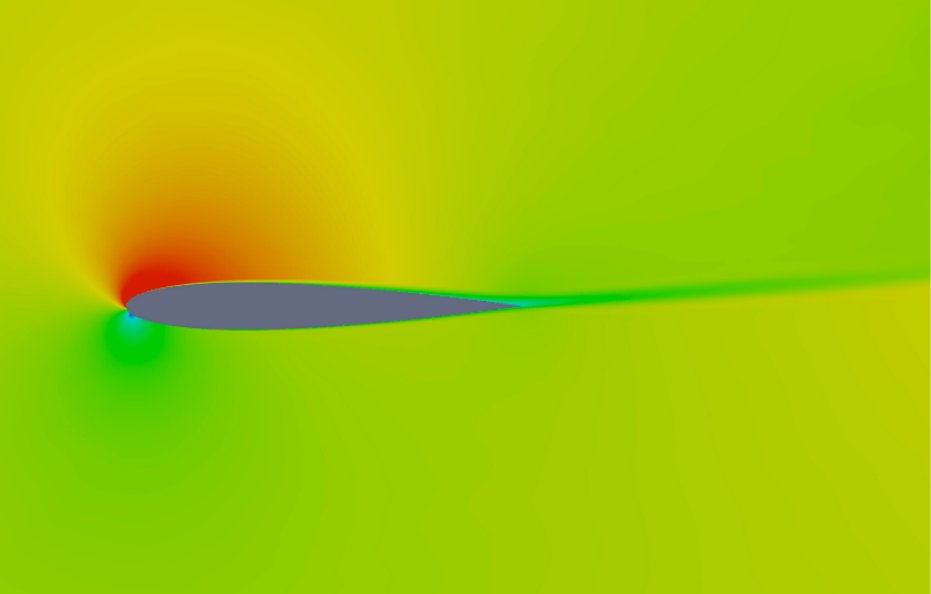}
  		\caption{Magnitude of Mean Velocity $\overline{{\bf{u}}}$}
	\end{subfigure}
		\begin{subfigure}{0.49\linewidth}
		\centering
  		\includegraphics[scale=0.85]{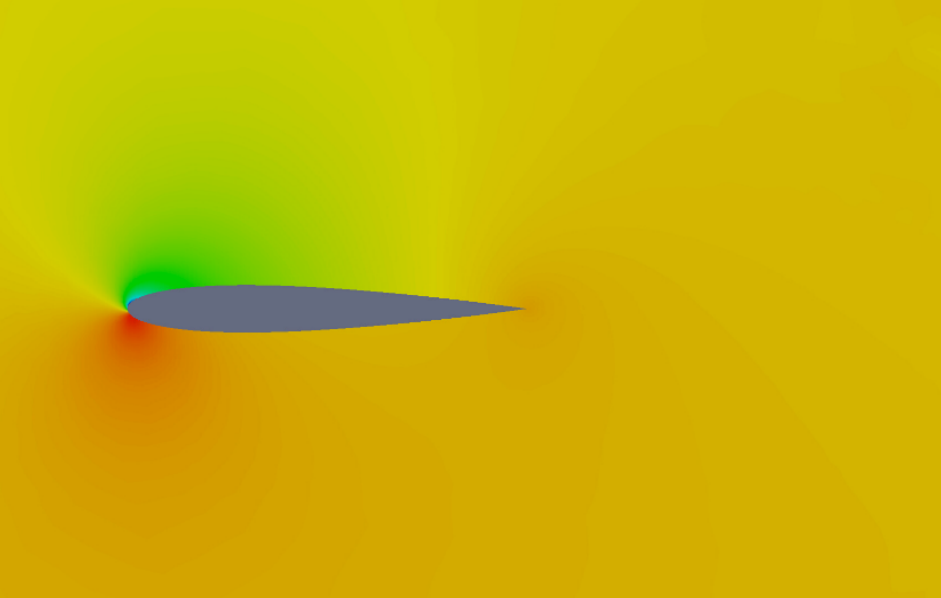}
  		\caption{Mean Pressure $\overline{p}$}
	\end{subfigure}
	\begin{subfigure}{1.0\linewidth}
		\centering
  		\includegraphics[scale=0.85]{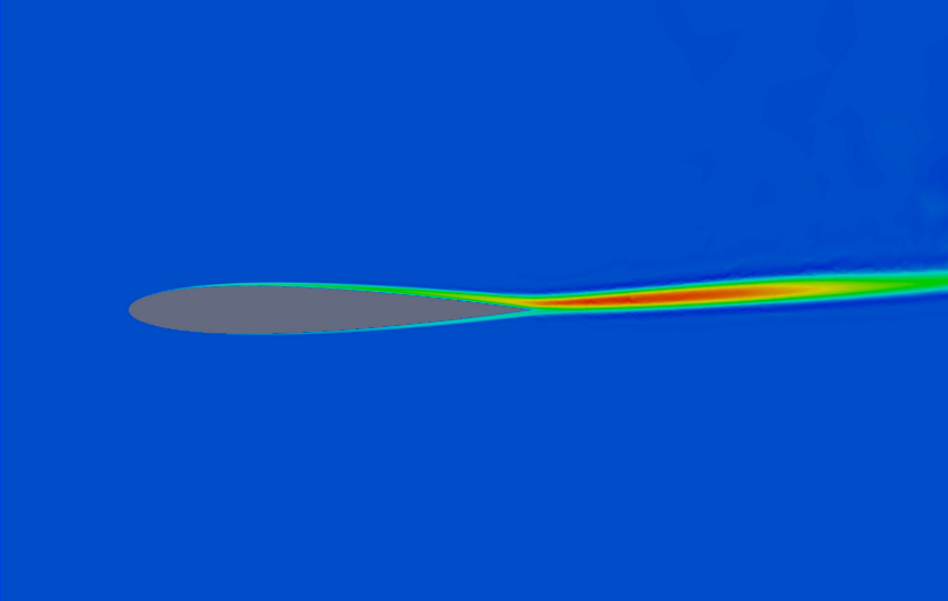}
  		\caption{Normalized Eddy Viscosity $\frac{\nu_T}{\nu}$}
	\end{subfigure}
	\caption{Computed velocity, pressure, and normalized eddy viscosity contours for flow over the NACA 0012 airfoil at $Re = 3,000,000$ and $\alpha=10^{\circ}$ using mean velocity, mean pressure, and working viscosity polynomial degrees of 2, 1, and 1 respectively and the mesh displayed in Fig. \ref{nacamesh}.}
	\label{airfoil10}
\end{figure}

\begin{figure}[t!]
    \centering
    \begin{subfigure}[b]{0.49\textwidth}
        \centering \includegraphics[scale  = 1]{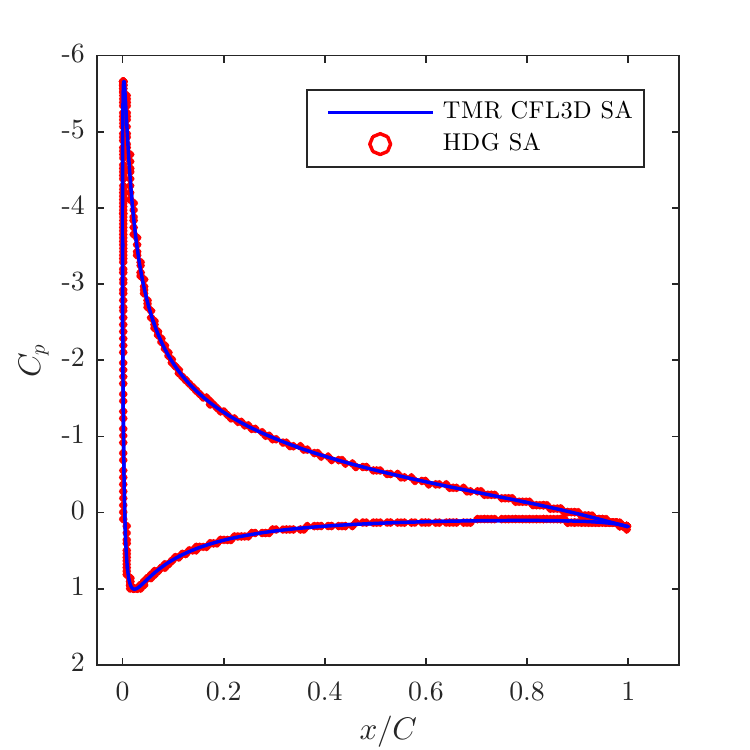}
        \caption{Surface Pressure Coefficient}
    \end{subfigure}
    \begin{subfigure}[b]{.49\textwidth}
        \centering \includegraphics[scale = 1]{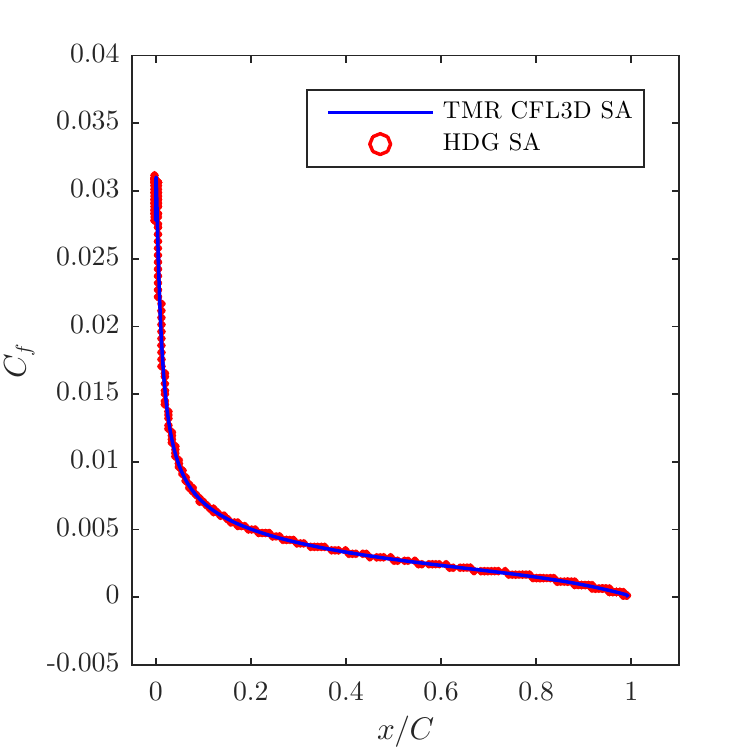}
        \caption{Surface Skin Friction Coefficient}
    \end{subfigure}
    \caption{Computed surface pressure and skin friction coefficients for flow over the NACA 0012 airfoil at $Re = 3,000,000$ and $\alpha=10^{\circ}$ using mean velocity, mean pressure, and working viscosity polynomial degrees of 2, 1, and 1 respectively and the mesh displayed in Fig. \ref{nacamesh}.  HDG results are plotted with finite volume results obtained using CFL3D and a mesh of 230,529 elements.}
        \label{cpcf10}
\end{figure}

As a final numerical experiment, we analyze turbulent flow over a NACA 0012 airfoil at a Reynolds number of $Re = 3,000,000$ based on the chord length and inflow velocity and an angle of attack of $\alpha = 10^{\circ}$.  This particular example has been explored in great depth in previous studies, and therefore there is an abundance of data to compare with for verification.  We specifically setup the problem to be similar to that described on the NASA turbulence modeling website \cite{nasa}.  The domain is chosen to be a sufficiently large box about the airfoil so that the boundary conditions do not affect the flow conditions near the airfoil.  The kinematic viscosity is set to $\nu = 3.33333e-7$, and the mean velocity field is set to be $\overline{{\bf{u}}} = (0.9848,0.1736)$ along the left (inflow) boundary.  Zero-traction boundary conditions are applied along the upper, lower, and right (outflow) boundaries, and homogeneous no-slip boundary conditions for the velocity field are applied along the airfoil surface itself.  For the working viscosity, four times the kinematic viscosity is prescribed as a Dirichlet inflow condition at the left boundary, a zero Dirichlet boundary condition is applied to the surface of the airfoil, and a zero traction Neumann condition is applied at the upper, lower, and right boundaries.\\
\\
To solve this flow problem, a computational mesh consisting of 61,022 elements and mean velocity, mean pressure, and working viscosity polynomial degrees of 2, 1, and 1, respectively, were employed.  A conservative grid stretching ratio of approximately 1.2 was used to generate the mesh, and special care was taken to ensure that the first point off the wall has a $y+$ value of approximately 1.  The computational mesh is displayed in Fig. \ref{nacamesh}, and attained steady state solutions for the mean velocity, mean pressure, and normalized eddy viscosity are displayed in Fig. \ref{airfoil10}.  The attained steady state solutions match well with steady solutions appearing in the literature \cite{burgess_high-order_2012}.  To further validate our results, we compare in Fig. \ref{cpcf10} the obtained surface pressure and skin friction coefficients on the airfoil with results publicly available on the NASA turbulence modeling website which were obtained using the NASA code CFL3D.  The CFL3D results were attained using a second-order finite volume method and a mesh of 230,529 elements.  The HDG and CFL3D results are visually indistinguishable, confirming the accuracy of the HDG results even though the HDG simulation only employed a quarter of the number of elements as the CFL3D simulation.

\section{Conclusions}
In this paper, we have introduced a new hybridized discontinuous Galerkin (HDG) method for the incompressible Reynolds Averaged Navier-Stokes (RANS) equations coupled with the Spalart-Allmaras one equation turbulence model, which extends the earlier work of Rhebergen and Wells in \cite{rhebergen_hybridizable_2017}.  We proved that our method is consistent, returns a point-wise divergence-free mean velocity field, conserves momentum globally, and is energy stable provided the turbulent eddy viscosity is non-negative.  We further demonstrated how static condensation could be employed to dramatically reduce the computational cost associated with each time step.  By turning to the method of manufactured solutions, we confirmed that our method yields optimal convergence rates provided the polynomial degree of the working viscosity is at most one order lower than that of the mean velocity field, and we further demonstrated the accuracy of our method using three benchmark problems: flow in a turbulent channel, flow over a backward facing step, and flow over a NACA 0012 airfoil at moderate Reynolds number and angle of attack.

There are several directions that we propose to explore in future work.  First of all, we will further study the influence of the polynomial degree of the working viscosity on the accuracy of the mean velocity and pressure fields, and we will also explore the potential of using different computational meshes for the flow field and the working viscosity.  Second, we plan to examine multilevel linear and nonlinear system solution strategies with an eye toward continually improving computational efficiency.  It should be noted that efficient multilevel solution strategies have been proposed for HDG methods in other application areas \cite{cockburn_multi, fabien_multigrid}.  Third, we plan to extend the method presented here to other turbulence models, including the $k$-$\varepsilon$, $k$-$\omega$, and SST models.  Finally, we plan to incorporate additional physics in future work such as chemical kinetics.

\section{Acknowledgements}
This material is based upon work supported by the Air Force Office of Scientific Research under Grant No. FA9550-14-1-0113.


\bibliography{main.bbl}

\end{document}